\documentclass[a4paper,11pt]{amsart}

\usepackage{a4wide}
\usepackage{amsmath, amsfonts, amscd, amssymb, amsthm}
\usepackage{tikz}
\usepackage[UKenglish]{babel}
\usetikzlibrary{matrix,arrows,calc}
\usepackage{xspace}

\newcommand{\thedocumentname}{Fibrations and log-symplectic structures}
\newcommand{\theauthor}{Gil R.\ Cavalcanti and Ralph L.\ Klaasse}
\usepackage{aliascnt}
\usepackage[pdftitle={\thedocumentname},pdfauthor={\theauthor},bookmarks,urlcolor=blue,colorlinks=true,linkcolor=black,citecolor=black]{hyperref}

\theoremstyle{plain}
\newtheorem{thm}{Theorem}[section]

\newtheorem*{thm2}{Theorem}
\newaliascnt{lem}{thm}
\newtheorem{lem}[lem]{Lemma}
\aliascntresetthe{lem}

\newaliascnt{prop}{thm}
\newtheorem{prop}[prop]{Proposition}
\aliascntresetthe{prop}

\newaliascnt{cor}{thm}
\newtheorem{cor}[cor]{Corollary}
\aliascntresetthe{cor}

\theoremstyle{definition}
\newaliascnt{defn}{thm}
\newtheorem{defn}[defn]{Definition}
\aliascntresetthe{defn}

\newaliascnt{rem}{thm}
\newtheorem{rem}[rem]{Remark}
\aliascntresetthe{rem}

\newaliascnt{exa}{thm}
\newtheorem{exa}[exa]{Example}
\aliascntresetthe{exa}

\newcommand{\N}{\mathbb{N}}
\newcommand{\Z}{\mathbb{Z}}
\newcommand{\R}{\mathbb{R}}
\newcommand{\C}{\mathbb{C}}
\newcommand{\im}{{\rm im}\,}
\newcommand{\bi}{\begin{itemize}}
\newcommand{\ei}{\end{itemize}}
\newcommand{\be}{\begin{equation*}}
\newcommand{\ee}{\end{equation*}}
\newcommand{\bp}[1][]{\begin{proof}[Proof#1.]}
\newcommand{\ep}{\end{proof}}
\newcommand{\wt}{\widetilde}
\newcommand{\wh}{\widehat}
\newcommand{\ol}{\overline}
\newcommand{\mc}{\mathcal}
\newcommand{\symp}{symplectic structure\xspace}

\newcommand{\blog}{log-\symp}
\newcommand{\Blog}{Log-\symp}
\newcommand{\bsymp}{$b$-\symp}

\newcommand{\fsymp}{folded-\symp}
\newcommand{\Fsymp}{Folded-\symp}
\newcommand{\nsymp}{near-\symp}
\newcommand{\Nsymp}{Near-\symp}
\newcommand{\acs}{almost-complex structure\xspace}

\newcommand{\bacs}{$b$-\acs}
\newcommand{\lf}{Lefschetz fibration\xspace}
\newcommand{\alf}{achiral \lf}
\newcommand{\Alf}{Achiral \lf}
\newcommand{\blf}{$b$-\lf}
\renewcommand{\b}[1][]{{^b}{#1}}
\newcommand{\btwo}[1][]{{^b}{#1}}
\renewcommand{\btwo}[1][]{#1}
\newcommand{\pt}[1][x]{\partial_{#1}}

\numberwithin{equation}{section}
\begin{document}

\author{Gil R.\ Cavalcanti}
\address{Department of Mathematics, Utrecht University, 3508 TA Utrecht, The Netherlands}
\email{g.r.cavalcanti@uu.nl}

\author{Ralph L.\ Klaasse}
\address{Department of Mathematics, Utrecht University, 3508 TA Utrecht, The Netherlands}
\email{r.l.klaasse@uu.nl}

\date{\today}
\thanks{This project was supported by a VIDI grant from NWO, the Dutch science foundation.}

\begin{abstract} \Blog{}s are Poisson structures $\pi$ on $X^{2n}$ for which $\bigwedge^n \pi$ vanishes transversally. By viewing them as symplectic forms in a Lie algebroid, the $b$-tangent bundle, we use symplectic techniques to obtain existence results for \blog{}s on total spaces of fibration-like maps. More precisely, we introduce the notion of a $b$-hyperfibration and show that they give rise to \blog{}s. Moreover, we link \blog{}s to \alf{}s and \fsymp{}s.
\end{abstract}

\title{Fibrations and log-symplectic structures}
\maketitle

\vspace{-2em}
\tableofcontents
\vspace{-3.5em}
\section{Introduction}
\label{sec:introduction}
A \blog{} is a Poisson structure $\pi$ on a $2n$-dimensional manifold $X$ for which $\bigwedge^n \pi$ is transverse to the zero section in $\bigwedge^{2n} TX$. \Blog{}s are Poisson structures that are in some sense as close as possible to being symplectic: if $\bigwedge^n \pi$ were nowhere vanishing, one could invert $\pi$ to obtain a nondegenerate closed two-form $\omega = \pi^{-1}$, i.e.\ a \symp on $X$. Transversality is a natural condition to impose on a Poisson structure and hence one is lead to study these objects. From this point of view, \blog{}s can be grouped with other symplectic-like structures whose singular behaviour is determined by a transversality condition, such as folded-symplectic and near-symplectic structures. The transversality condition implies that instead of $\pi$ being invertible, it drops rank on -- at most -- a codimension-one submanifold $Z_\pi = (\bigwedge^n \pi)^{-1}(0)$, and there does so generically. The existing literature on \blog{}s includes \cite{Cavalcanti13, FrejlichMartinezTorresMiranda15, GMP13, MarcutOsorno14, MarcutOsorno14two, Radko02}, where they sometimes go under the name of $b$-Poisson or $b$-\symp{}s.

In this paper we take the following viewpoint: out of the singular locus $Z_\pi$ one defines a Lie algebroid, the $b$-tangent bundle $\b{TX}$. One then notes that \blog{}s with given $Z_\pi$ are nothing more than \symp{}s for $\b{TX}$ (see Section \ref{sec:blogsbtangent} for further details). Through this use of $b$-geometry one sees that \blog{}s are close enough to being symplectic that one can use symplectic techniques to study them.

The main goal of this paper is to use the language of $b$-geometry to extend results from symplectic geometry relating the existence of fibration-like maps to \blog{}s. We first give a brief overview of what has been done in this direction.

In \cite{Thurston76}, see also \cite{McDuffSalamon98}, Thurston showed how to equip symplectic fiber bundles with \symp{}s. Gompf \cite{GompfStipsicz99} then showed that Lefschetz fibrations lead to \symp{}s in dimension four, adapting Thurston's methods. Conversely, Donaldson \cite{Donaldson99} proved using approximately holomorphic methods that symplectic forms lead to Lefschetz pencils. Further, Gompf \cite{Gompf04two, Gompf04} introduced hyperpencils in all dimensions and showed they admit \symp{}s, aiming to give a topological characterization of symplectic manifolds. Indeed, having established this correspondence, the study of Lefschetz pencils can shed light onto symplectic geometry. One can then branch out to other symplectic-like structures, and connect them to suitably generalized Lefschetz-type fibrations. This has been done for near-symplectic and folded-symplectic structures.
\subsection*{Near-symplectic geometry}
Recall that a \nsymp on a compact oriented four-manifold $X$ is a closed two-form $\omega$ such that there exists a Riemannian metric making $\omega$ harmonic, self-dual, and transverse to the zero section in $\bigwedge^{2,+} T^* X$. Such forms vanish on circles, and are symplectic whenever they do not vanish. \Nsymp{}s exist on all compact oriented four-manifolds with $b_2^+ > 0$, as one can prove by studying generic metrics \cite{Honda04,LeBrun97}. Auroux--Donaldson--Katzarkov \cite{AurouxDonaldsonKatzarkov05} introduced broken Lefschetz pencils to study \nsymp{}s, and established the same correspondence as between Lefschetz pencils and \symp{}s.

Instead of allowing for broken-type singularities, one can drop the requirement that the charts around Lefschetz singularities are orientation preserving, obtaining so-called achiral Lefschetz singularities. Correspondingly, there is the notion of an \alf{} and its broken variant. Gay--Kirby \cite{GayKirby07} proved that all closed oriented four-manifolds admit broken \alf{}s. Soon thereafter Baykur \cite{Baykur08}, Lekili \cite{Lekili09} and Akbulut--Karakurt \cite{AkbulutKarakurt08} showed the achirality is unnecessary, i.e.\ that one can always find broken \lf{}s over the two-sphere, and broken Lefschetz pencils if in addition $b_2^+ > 0$. See also \cite{Baykur09} for more on broken \lf{}s and \nsymp{}s.

Returning to \alf{}s, Etnyre--Fuller \cite{EtnyreFuller06} showed that any oriented closed four-manifold admits an \alf over the two-sphere, after performing surgery on a circle. Moreover, they showed that any \alf with a section gives rise to a \nsymp.
\subsection*{Folded-symplectic geometry}
Extending the symplectic condition in another direction one obtains \fsymp{}s. These are closed two-forms $\omega$ on a $2n$-dimensional manifold $X$ for which $\omega^n$ is transverse to the zero section in $\bigwedge^{2n} T^*X$, and such that the restriction of $\omega^{n-1}$ to the zero set $Z_\omega$ of $\omega^n$ has maximal rank. Here, as for \blog{}s, $Z_\omega$ is a hypersurface of $X$. Cannas da Silva showed that \fsymp{}s exist on all compact oriented four-manifolds \cite{CannasDaSilva10}, as their existence is equivalent to having a stable \acs. This was proven using a version of the $h$-principle for folding maps. In \cite{Baykur06}, Baykur gave a topological proof for their existence, and showed they arise out of \alf{}s.
\begin{rem}\label{rem:homessential} In all of the above instances, having some form of homological essentialness of the fibers is required, as the symplectic-like structure will restrict to a symplectic (volume) form on the fibers $F$ of the Lefschetz-type fibration. This is guaranteed for Lefschetz fibrations with fiber genus different from one (as the kernel of the map defines a line bundle $L$, whose Chern class satisfies $\langle c_1(L)|_F, [F] \rangle = \chi(F) \neq 0$). Furthermore, it is guaranteed whenever the fibration-like map admits a section. However, in the general case, one must impose this condition by hand. There are in some sense two ways of doing this, one of which is by asking that $[F] \neq 0$ and ensuring the fibers are connected. Alternatively, by Poincar\'e duality, one can demand the existence of a closed two-form which pairs positively with every fiber component. When dealing with pencils this is usually not required, as one can then blow-up the base locus to obtain sections of the resulting fibration.
\end{rem}
\subsection*{Log-symplectic geometry}
As was mentioned before, we study \blog{}s by viewing them as symplectic forms in a Lie algebroid, called the $b$-tangent bundle, and use this viewpoint to apply constructions from symplectic geometry. In doing so we obtain existence results for \blog{}s on total spaces of fibration-like maps. We say a pair $(X,Z_X)$ consisting of a manifold and a hypersurface admits a \blog $\pi$ if $Z_\pi = Z_X$, and similarly for \fsymp{}s. As a warm-up, we prove the $b$-analogue of Thurston's result for fibrations with two-dimensional fibers as \autoref{thm:bthurstonfibration}, which implies the following (see \autoref{cor:fibrationblogs}).
\begin{thm2} Let $f:\, X^{2n} \to Y^{2n-2}$ be a fibration between compact connected manifolds. Assume that $Y$ admits a \blog $\pi$ and that the generic fiber $F$ is orientable and $[F] \neq 0 \in H_2(X;\R)$. Then $(X,Z_X)$ admits a \blog, where $Z_X = f^{-1}(Z_{\pi})$.
\end{thm2}
Further, after introducing the notion of a $b$-hyperfibration (the $b$-analogue of a hyperpencil with empty base locus \cite{Gompf04}) in Section \ref{sec:bhyperfibration}, we show the following as \autoref{thm:bhyperfibration}.
\begin{thm2}\label{thm:introbhyperfibration} Let $f:\, (X,Z_X) \to (Y, Z_Y, \btwo{\omega}_Y)$ be a $b$-hyperfibration between compact connected $b$-oriented $b$-manifolds. Assume that there exists a finite collection $S$ of sections of $f$ interescting all fiber components non-negatively and for each fiber component at least one section in $S$ intersecting positively. Then $(X,Z_X)$ admits a \blog.
\end{thm2}
Our methods also give an alternative proof of a result in \cite{Cavalcanti13} that \alf{}s on compact four-manifolds lead to \blog{}s, which we show as \autoref{thm:alflogsymp}.
\begin{thm2}\label{thm:introalflogsymp} Let $f:\, X^4 \to \Sigma^2$ be an \alf between compact connected manifolds. Assume that the generic fiber $F$ is orientable and $[F] \neq 0 \in H_2(X;\R)$. Then $X$ admits a \blog.
\end{thm2}
It was mentioned before that one can obtain a \fsymp out of an \alf, but from our point of view this is due to the theorem above and the fact that any \blog induces a \fsymp. This result is known in the Poisson community, and we include its proof here (see \autoref{thm:blogfolded}) for completeness. Further, we show when one can obtain a \blog out of a \fsymp as \autoref{thm:foldedblog}, establishing a converse.
\begin{thm2}\label{thm:introblogfolded} Let $X^{2n}$ be a compact manifold and $Z_X$ a hypersurface. Then $(X,Z_X)$ admits a \blog $\pi$ if and only if $(X,Z_X)$ admits a \fsymp $\omega$ and there exists a closed one-form $\theta \in \Omega^1(Z_X)$ satisfying $\theta \wedge \omega^{n-1}|_{Z_X} \neq 0$. Moreover, $\pi^{-1} = \omega$ outside a neighbourhood of $Z_X$.
\end{thm2}
While every oriented four-manifold admits \fsymp{}s \cite{CannasDaSilva10}, not every compact oriented four-manifold admits a \blog, as there are cohomological obstructions to the existence of the latter, see \cite{Cavalcanti13,MarcutOsorno14two} or Section \ref{sec:blogsbtangent}.

All manifolds in this paper will be compact and without boundary, unless specifically stated otherwise. Note however that they are not necessarily oriented nor orientable. Throughout we will identify de Rham cohomology and singular cohomology with $\R$-coefficients.
\subsection*{Organization of the paper}
In Section \ref{sec:blogsbtangent} we discuss the definition of \blog{}s and describe how they can be viewed as symplectic forms for the $b$-tangent bundle. We further discuss the existence of \blog{}s on surfaces. In Section \ref{sec:blogconstructions} we prove \autoref{thm:bthurstontrick}, which is our main tool for constructing \blog{}s, and allows us to prove the $b$-analogue of Thurston's result, \autoref{thm:bthurstonfibration}. In Section \ref{sec:achirallefschetz} we use this tool to prove \autoref{thm:alflogsymp} that four-dimensional \alf{}s give rise to \blog{}s. In Section \ref{sec:bhyperfibration} we define the notion of a $b$-hyperfibration and prove our main theorem,  \autoref{thm:bhyperfibration}, that $b$-hyperfibrations lead to \blog{}s. Finally, Section \ref{sec:blogsandfsymps} discusses the relation between log-symplectic and \fsymp{}s.
%
%
\section{\Blog{}s and the \texorpdfstring{$b$}{b}-tangent bundle}
\label{sec:blogsbtangent}
In this section we discuss the notion of a \blog{} and the $b$-geometry language used to study these structures. We defer a more comprehensive account of the interplay of $b$-geometry with \blog{}s to \cite{Klaasse16} and instead only develop what is required for our purposes.
\begin{defn} A \emph{Poisson bivector} $\pi$ on a compact manifold $X$ is a section $\pi \in \Gamma(\bigwedge^2 TX)$ satisfying $[\pi, \pi] = 0$, where $[\cdot,\cdot]$ is the Schouten-Nijenhuis bracket of multivector fields.
\end{defn}
\begin{defn}\label{defn:blog} A \emph{\blog} on a compact $2n$-dimensional manifold $X$ is a Poisson bivector $\pi$ such that $\bigwedge^n \pi$ is transverse to the zero section in $\bigwedge^{2n} TX$. The set $Z_\pi = (\bigwedge^n \pi)^{-1}(0)$ is called the \emph{singular locus} of $\pi$, and a \blog is called \emph{bona fide} if $Z_\pi \neq \emptyset$.
\end{defn}
The singular locus of a \blog is a codimension-one smooth submanifold, or \emph{hypersurface}, of $X$. It may be empty and need not be connected. If a given Poisson bivector $\pi$ would be of full rank, we could invert it to obtain a symplectic structure $\omega = \pi^{-1}$ on $X$. One could thus say that \blog{}s are ``generically generic'' Poisson structures, as they drop rank on at most a hypersurface in $X$ and there do so generically. We call $X \setminus Z_\pi$ the \emph{symplectic locus} of $\pi$. We are interested in equipping manifolds with \blog{}s and proceed by viewing them as symplectic structures in a specific Lie algebroid. This point of view has also been adopted by \cite{FrejlichMartinezTorresMiranda15, GMP13, MarcutOsorno14} and others.

Let $(X,Z_X)$ be a \emph{pair}, i.e.\ a manifold $X$ together with a hypersurface $Z_X \subset X$. We say the pair $(X,Z_X)$ \emph{admits a \blog} if there exists a \blog $\pi$ on $X$ such that $Z_\pi = Z_X$. Let $\mc{V}_b(X) \subset \Gamma(TX)$ be the set of vector fields on $X$ which are tangent to $Z_X$. Then $\mc{V}_b(X)$ defines a locally free sheaf, hence leads to a vector bundle by the Serre-Swan theorem. Furthermore, one notes that $\mc{V}_b$ forms a Lie subalgebra inside $\Gamma(TX)$.
\begin{defn} Let $(X,Z_X)$ be a pair. The \emph{$b$-tangent bundle} $\b{TX} \to X$ is the vector bundle on $X$ with $\Gamma(\b{TX}) = \mc{V}_b(X)$. A \emph{$b$-manifold} is a pair $(X,Z_X)$ equipped with the bundle $\b{TX}$.
\end{defn}
We will not distinguish between pairs and $b$-manifolds. The $b$-tangent bundle is an example of a \emph{Lie algebroid}, i.e.\ a vector bundle $L \to X$ together with a map $\rho_X:\, L \to TX$ called the \emph{anchor}, such that $\Gamma(L)$ is a Lie algebra, $\rho_X$ is a map of Lie algebras when regarded as a map between spaces of sections, and $\rho_X$ satisfies the Leibniz rule $[v, f w] = f [v, w] + (\rho_X(v) f) w$ for all $v, w \in \Gamma(L)$ and $f \in C^\infty(X)$. In this case the anchor is the natural inclusion, which is an isomorphism away from $Z_X$. The $b$-tangent bundle in the case when $X$ is a manifold with boundary and $Z_X = \partial X$ has been extensively studied by Melrose and others \cite{Melrose93}.
\begin{defn} Given two $b$-manifolds $(X,Z_X)$ and $(Y,Z_Y)$, a \emph{$b$-map} is a map $f:\, X \to Y$ such that $f^{-1}(Z_Y) = Z_X$ and $f$ is transverse to $Z_Y$. We write $f:\, (X,Z_X) \to (Y,Z_Y)$.
\end{defn}
In other words, for a $b$-map we have $df_x(T_x X) + T_y Z_Y = T_y Y$ for all $y \in Z_Y$, where $x \in f^{-1}(y)$. Note that $Z_X$ is uniquely determined by $Z_Y$ and the requirement that $f$ is a $b$-map. One checks that there is a category with objects being $b$-manifolds and morphisms being $b$-maps between them. Given a $b$-map $f:\, (X,Z_X) \to (Y,Z_Y)$, its level sets are either contained in $Z_X$ or are disjoint from it. The anchor $\rho_X$ gives rise to a well-defined line bundle $\mathbb{L}_X := \ker \rho_X$ over $Z_X$ which is always trivial. Indeed, $\mathbb{L}_X$ is canonically trivialized by a nonvanishing section called the \emph{normal $b$-vector field} \cite[Proposition 4]{GMP13}, locally given by $x \frac{\partial}{\partial x}$ for $x$ a local defining function for $Z_X$.
\begin{defn} A \emph{$b$-fibration} is a surjective $b$-submersive $b$-map $f:\, (X,Z_X) \to (Y,Z_Y)$.
\end{defn}
\begin{rem}\label{rem:fibrintobfibr} It follows immediately from the definition that given a fibration $f:\, X^{2n} \to Y^{2n-2}$, one can turn it into a $b$-fibration by choosing a hypersurface $Z_Y \subset Y$ and considering the $b$-map $f:\, (X,Z_X) \to (Y,Z_Y)$, where $Z_X = f^{-1}(Z_Y)$.
\end{rem}
There are the associated notions of Lie algebroid $k$-forms which we call \emph{$b$-$k$-forms} and denote their space of sections by $\b{\Omega}^k(X)$, and de Rham differential $\b{d}:\, \b{\Omega}^k(X) \to \b{\Omega}^{k+1}(X)$ giving the Lie algebroid de Rham cohomology $\b{H}^k(X)$. Note that, given a $b$-map $f:\, (X,Z_X) \to (Y,Z_Y)$ and a $b$-form $\btwo{\omega}_Y \in \b{\Omega}^k(Y)$, the pullback $f^* \btwo{\omega}_Y$ is a well-defined element of $\b{\Omega}^k(X)$. We can view forms on $X$ as $b$-forms using pullback via the anchor, $\rho_X^*:\, \Omega^k(X) \to \b{\Omega}^k(X)$.
\begin{rem} For the sake of readability, we will often not typographically distinguish between $k$-forms and $b$-$k$-forms, as well as write merely $d$ instead of $\b{d}$ and say closed instead of $b$-closed. This should not cause any confusion, but note for example that for forms viewed as $b$-forms, the notions of exactness and $b$-exactness do not agree.
\end{rem}
A \emph{$b$-orientation} for $(X,Z_X)$ is an orientation for the bundle $\b{TX}$. A \emph{\bacs} is a vector bundle complex structure for $\b{TX}$, which induces a $b$-orientation. Note that a $b$-manifold $(X,Z_X)$ may be $b$-orientable yet not be orientable, and vice versa.
\begin{defn} A \emph{$b$-symplectic form} is a closed, nondegenerate $b$-two-form $\btwo{\omega} \in \b{\Omega}^2(X)$.
\end{defn}
As in the symplectic case, a $b$-symplectic form induces a $b$-orientation. We will only consider $b$-orientable $(X,Z_X)$, and \bacs{}s and $b$-symplectic forms inducing the same $b$-orientation. The reason for introducing the $b$-tangent bundle and $b$-symplectic forms is the following. Viewing a $b$-two-form $\btwo{\omega} \in \b{\Omega}^2(X)$ as a map $\btwo{\omega}:\, \b{TX} \to \b{T^*X}$, if $\btwo{\omega}$ is of maximal rank, this map can be inverted to a map $\b{\pi} := \btwo{\omega}^{-1}:\, \b{T^*X} \to \b{TX}$ given by some $b$-bivector $\b{\pi} \in \Gamma(\bigwedge^2 \b{TX})$. Using the anchor one obtains $\rho_X \circ \b{\pi}:\, \b{T^*X} \to TX$, specifying a bivector $\pi_{\btwo{\omega}} := \rho_X(\b{\pi}) \in \Gamma(\bigwedge^2 TX)$, which is called the \emph{dual bivector} to $\btwo{\omega}$. On the other hand, given a \blog $\pi$ on $X$, one sees that $\pi:\, T^*X \to \b{TX}$, and there is a unique $b$-bivector $\b{\pi}:\, \b{T^*X} \to \b{TX}$ such that $\pi = \b{\pi} \circ \rho_X^*$. The dual $b$-two-form $\btwo{\omega} = \b{\pi}^{-1}$ is then a $b$-symplectic form. This is illustrated by the following two diagrams.

\begin{center}
\begin{tikzpicture}
\begin{scope}[xshift=-160pt]
  \matrix (m) [matrix of math nodes,row sep=3em,column sep=4em,minimum width=2em]
  {
     \b{TX} & \b{T^*X}\\
     TX & T^*X\\};
   \path[-stealth,transform canvas={yshift=2pt}]
   (m-1-1) edge node [above] {$\btwo{\omega}$} (m-1-2);
   \path[-stealth,transform canvas={yshift=-2pt}]
   (m-1-2) edge node [below] {$\b{\pi}$} (m-1-1);
   \path[-stealth]
    (m-1-1) edge node [left] {$\rho_X$} (m-2-1)
    (m-2-2) edge node [below] {$\pi_{\btwo{\omega}}$} (m-2-1)
                 edge node [right] {$\rho_X^*$} (m-1-2);
\end{scope}
\begin{scope}
  \matrix (m) [matrix of math nodes,row sep=3em,column sep=4em,minimum width=2em]
  {
     T^* X & \b{TX}\\
     \b{T^*X}\\};
  \path[-stealth]
    (m-1-1) edge node [above] {$\pi$} (m-1-2)
                 edge node [left] {$\rho_X^*$} (m-2-1)
    (m-2-1) edge node [below right] {$\b{\pi}$} (m-1-2);
\end{scope}
\end{tikzpicture}
\end{center}

These processes are inverse to each other, and we have the following key result.
\begin{prop}[{\cite[Proposition 20]{GMP13}}]\label{prop:blogbsymp} Given a manifold $X$, a $b$-two-form $\btwo{\omega}$ on the $b$-manifold $(X,Z_X)$ for some hypersurface $Z_X$ is $b$-symplectic if and only if its dual bivector $\pi = \pi_{\btwo{\omega}}$ is a \blog with $Z_{\pi} = Z_X$. In particular, if $\pi$ is bona fide log-symplectic, then $Z_X \neq \emptyset$.
\end{prop}
We thus see that the pair $(X,Z_X)$ admits a \blog if and only if the $b$-manifold $(X,Z_X)$ admits a \bsymp. The $b$-two-form associated to a \blog $\pi$ by the above procedure can be viewed as a symplectic form which has a logarithmic singularity at $Z_X$. Indeed, the local Darboux model for a $b$-symplectic form $\btwo{\omega}$ on a $2n$-dimensional manifold $X$ is given by $\btwo{\omega} = d \log x_1 \wedge d x_2 + \dots + d x_{2n-1} \wedge d x_{2n}$, using coordinates $x_i$ in a neighbourhood $U$ such that $Z_X \cap U = \{x_1 = 0\}$. Its dual bivector is given locally by $\pi = x_1 \pt[x_1] \wedge \pt[x_2] + \dots + \pt[x_{2n-1}] \wedge \pt[x_{2n}]$, where $\pt[x_i] = \partial / \partial x_i$.
\begin{rem} As $\b{TX}$ is isomorphic to $TX$ away from $Z_X$ by the anchor, a $b$-oriented $b$-manifold $(X,Z_X)$ obtains an orientation in the usual sense away from $Z_X$. Given a $b$-symplectic manifold $(X,Z_X,\btwo{\omega})$ with dual bivector $\pi$ and the induced $b$-orientation, this orientation on $X \setminus Z_X$ can never come from an already existing orientation on $X$ when $Z_X \neq \emptyset$. Either $X$ is non-orientable, or if $s \in \Gamma(\bigwedge^{2n} TX)$ is a nonvanishing section orienting $X$, we have $\bigwedge^n \pi = h s$ for some smooth function $h$ vanishing precisely at $Z_X$. From the local description of $\pi$ above one sees that $h$ must change sign so that $Z_X$ separates $X$ into two connected components according to the sign of $h$.
\end{rem}
The following result extracted from \cite{GMP13} describes the inverse of a \blog around its singular locus, see also \cite[Theorem 3.2]{Cavalcanti13}.
\begin{prop}\label{prop:bloglocalform} Let $(X^{2n},Z_X,\pi)$ be a compact log-symplectic manifold. Then around any connected component $Z$ of $Z_X$, the two-form $\pi^{-1}$ is equivalent to $d \log |x| \wedge \theta + \sigma$ in a neighbourhood of the zero section of the normal bundle $NZ$ of $Z$, where $|x|$ is the distance to the zero section with respect to a fixed metric on $NZ$. Here $\theta$ and $\sigma$ are closed one- and two-forms on $Z$ satisfying $\theta \wedge \sigma^{n-1} \neq 0$.
\end{prop}
The above proposition thus also states that a \blog induces a \emph{cosymplectic structure} on its singular locus. Moreover, let us mention that there are cohomological obstructions for a compact manifold to admit a \blog.
\begin{thm}[\cite{MarcutOsorno14two}]\label{thm:blogobstr} Let $X^{2n}$ be a compact log-symplectic manifold. Then there exists a class $a \in H^2(X;\R)$ such that $a^{n-1} \neq 0$.
\end{thm}
\begin{thm}[\cite{Cavalcanti13}]\label{thm:bonafideobstr} Let $X^{2n}$ be a compact oriented bona fide log-symplectic manifold. Then there exists a nonzero class $b \in H^2(X;\R)$ such that $b^2 = 0$. Moreover, if $n > 1$ then $b_2(X) \geq 2$.
\end{thm}
We now discuss the extent to which there is a difference between $df$-critical points and $\b{df}$-critical points, given a $b$-map $f:\, (X,Z_X) \to (Y,Z_Y)$.
\begin{prop}\label{prop:kerbdf} Let $f:\, (X,Z_X) \to (Y,Z_Y)$ be a $b$-map. Then $\rho_{X,x}:\, \ker \b{df}_x \to \ker df_x$ is an isomorphism for all $x \in X$.
\end{prop}
Consequently we can unambiguously speak of a critical point of $f$, without specifying whether we mean with respect to $\b{df}$ or $df$. The essential ingredients are contained in the following lemma. Two vector spaces $V, V_1$ will be called a \emph{pair} if $V_1 \subset V$ is a subspace. A linear map $f:\, V \to W$ is a map between pairs $(V,V_1)$ and $(W,W_1)$ if $f(V_1) \subset W_1$.
\begin{lem}\label{lem:blinalg} Let $F:\, (V,V_1) \to (W,W_1)$ be a linear map between pairs such that under the projection maps ${\rm pr}_V:\, V \mapsto V / V_1$ and ${\rm pr}_W:\, W \mapsto W / W_1$, $F$ descends to an isomorphism $\overline{F}:\, V / V_1 \to W / W_1$. Assume that there are vector spaces $\b{V}$, $\b{W}$ and maps $\rho_V:\, \b{V} \to V$, $\rho_W:\, \b{W} \to W$, $\b{F}:\, \b{V} \to \b{W}$ so that $F \circ \rho_V = \rho_W \circ \b{F}$. Assume that $\im \rho_V = V_1$, $\im \rho_W = W_1$ and $\b{F}:\, \ker \rho_V \to \ker \rho_W$ is an isomorphism. Then $\rho_V: \, \ker \b{F} \to \ker F$ is an isomorphism.
\end{lem}
The situation is summarized by the following diagram, in which the rows are exact. The two vertical maps on the far left and right are assumed to be isomorphisms, while the conclusion of the lemma is that the top horizontal map is an isomorphism.
\begin{center}
\begin{tikzpicture}
\matrix (m) [matrix of math nodes, row sep=2.5em, column sep=2.5em,text height=1.5ex, text depth=0.25ex]
{
 &                      & \ker \b{F} & \ker F  &              &   \\
	0 & \ker \rho_V   & \b{V}        & V         & V/V_1   & 0 \\
	0 & \ker \rho_W  & \b{W}       & W        & W/W_1 & 0 \\};
\path[-stealth]
(m-1-3) edge node [above] {$\rho_V$} node [below] {$\cong$} (m-1-4)
edge [right hook-latex] (m-2-3)
(m-1-4) edge [right hook-latex] (m-2-4)
(m-2-1) edge (m-2-2)
(m-2-2) edge [right hook-latex] (m-2-3)
(m-2-2) edge node [left] {$\cong$} node [right] {$\b{F}$} (m-3-2)
(m-2-3) edge node [right] {$\b{F}$} (m-3-3)
edge node [above] {$\rho_V$} (m-2-4)
(m-2-4) edge node [right] {$F$} (m-3-4)
edge node [above] {${\rm pr}_V$} (m-2-5)
(m-2-5) edge (m-2-6)
(m-2-5) edge node [right] {$\overline{F}$} node [left] {$\cong$} (m-3-5)
(m-3-1) edge (m-3-2)
(m-3-2) edge [right hook-latex] (m-3-3)
(m-3-3) edge node [above] {$\rho_W$} (m-3-4)
(m-3-4) edge node [above] {${\rm pr}_W$} (m-3-5)
(m-3-5) edge (m-3-6);
\end{tikzpicture}
\end{center}
\bp If $\b{v} \in \ker \b{F}$, then $F \rho_V (\b{v}) = \rho_W \b{F} (\b{v}) = 0$, so that $\rho_V: \, \ker \b{F} \to \ker F$. Given $v \in \ker F$ we have ${\rm pr}_W F(v) = {\rm pr}_W(0) = 0$, so by assumption ${\rm pr}_V(v) = 0$, hence $v \in V_1$. As $v \in V_1$, there exists a $\b{v}_0 \in \b{V}$ such that $\rho_V(\b{v}_0) = v$, and $\rho_V^{-1}(v) = \b{v}_0 + \ker \rho_V$. Consider a vector $\b{v} = \b{v}_0 + k$ for $k \in \ker \rho_V$, so that $\rho_V (\b{v}) = v$. Then $\rho_W \b{F}(\b{v}) = F \rho_V(\b{v}) = F(v) = 0$, so that $\b{F}(\b{v}) \in \ker \rho_W$. But then there exists a unique $k \in \ker \rho_V$ such that $\b{F}(k) = - \b{F}(\b{v}_0)$, for which $\b{F}(\b{v)} = 0$. Hence there exists a unique $\b{v} \in \ker \b{F} \cap \rho_V^{-1}(v)$. We conclude that $\rho_V:\, \ker \b{F} \to \ker F$ is an isomorphism.
\ep
\bp[ of \autoref{prop:kerbdf}] Let $x \in X$ be given and denote $y = f(x)$. If $x \in X \setminus Z_X$ the statement follows as $\rho_X$ is an isomorphism away from $Z_X$. If $x \in Z_X$, we wish to apply \autoref{lem:blinalg} to the situation where $V = T_x X$, $V_1 = T_x Z_X$, $W = T_{y} Y$, $W_1 = T_{y} Z_Y$, $\b{V} = \b{T_x X}$, $\b{W} = \b{T_{y} Y}$, $F = df_x$, $\b{F} = \b{df}_x$, $\rho_V = \rho_{X,x}$ and $\rho_W = \rho_{Y,y}$. As $f$ is a $b$-map, we have $f^{-1}(Z_Y) = Z_X$ so that $df_x(T_x {Z_X}) \subset T_{y} Z_Y$ making $F$ a map between pairs. Splitting $T_x X = N_x Z_X \oplus T_x Z_X$ pointwise we see that $T_x X / T_x Z_X \cong N_x Z_X$ and $df_x(T_x X) + T_y Z_Y = df_x(N_x Z_X) + T_y Z_Y$. Because $f$ is a $b$-map we have $df_x(T_x X) + T_{y} Z_Y = T_{y} Y$. Similarly splitting $T_y Y = N_y Z_Y \oplus T_y Z_Y$ we conclude that $df_x(N_x Z_X) + T_y Z_Y = N_y Z_Y \oplus T_y Z_Y$ and we see that $df_x$ restricts to an isomorphism from $N_x Z_X \cong V / V_1$ to $N_y Z_Y \cong W / W_1$. By construction $df \circ \rho_X = \rho_Y \circ \b{df}$, hence also pointwise at $x$ and $y$, and furthermore by definition of the $b$-tangent bundles we have $\im \rho_{X,x} = T_x Z_X$ and $\im \rho_{Y,y} = T_{y} Z_Y$. Note that $f^* \mathbb{L}_Y \cong \mathbb{L}_X$, because local defining functions for $Z_Y$ pull back to local defining functions for $Z_X$. Hence $\b{df}_x:\, \mathbb{L}_{X,x} \to \mathbb{L}_{Y,y}$ is an isomorphism. Finally, $\ker \rho_{X,x} = \mathbb{L}_{X,x}$ and $\ker \rho_{Y,y} = \mathbb{L}_{Y,y}$, so that $\b{df}_x:\, \ker \rho_{X,x} \to \ker \rho_{Y,y}$ is an isomorphism. By \autoref{lem:blinalg} we conclude that $\rho_{X,x}:\, \ker \b{df}_x \to \ker df_x$ is an isomorphism.
\ep

The situation is summarized by the following diagram with exact rows.
\begin{center}
\begin{tikzpicture}
\matrix (m) [matrix of math nodes, row sep=2.5em, column sep=2.5em,text height=1.5ex, text depth=0.25ex]
{	   &                               & \ker \b{df}_x    & \ker df_x     &                             &   \\
	0 & \mathbb{L}_{X,x}   & \b{T_x X}        & T_x X         & T_x X / T_x Z_X   & 0 \\
	0 & \mathbb{L}_{Y,y}   & \b{T_y Y}        & T_y Y          & T_y Y / T_y Z_Y   & 0 \\};
\path[-stealth]
(m-1-3) edge node [above] {$\rho_{X,x}$} node [below] {$\cong$} (m-1-4)
edge [right hook-latex] (m-2-3)
(m-1-4) edge [right hook-latex] (m-2-4)
(m-2-1) edge (m-2-2)
(m-2-2) edge [right hook-latex] (m-2-3)
(m-2-2) edge node [left] {$\cong$} node [right] {$\b{df}_x$} (m-3-2)
(m-2-3) edge node [right] {$\b{df}_x$} (m-3-3)
edge node [above] {$\rho_{X,x}$} (m-2-4)
(m-2-4) edge node [right] {$df_x$} (m-3-4)
edge node [above] {${\rm pr}_x$} (m-2-5)
(m-2-5) edge[anchor = south] (m-2-6)
(m-2-5) edge node [right] {$\overline{df}_x$} node [left] {$\cong$} (m-3-5)
(m-3-1) edge (m-3-2)
(m-3-2) edge [right hook-latex] (m-3-3)
(m-3-3) edge node [above] {$\rho_{Y,y}$} (m-3-4)
(m-3-4) edge node [above] {${\rm pr}_y$} (m-3-5)
(m-3-5) edge (m-3-6);
\end{tikzpicture}
\end{center}
\begin{rem} The statement that $df_x:\, N_x Z_X \to N_y Z_Y$ is an isomorphism for $y = f(x) \in Z_Y$ can be colloquially phrased as follows. As $f$ is transverse to $Z_Y$, the normal direction to $Z_Y$ at points in $Z_Y$ must be contained in the image of $TX$ under $df$. The fact that $f^{-1}(Z_Y) = Z_X$ then implies that it must in fact be obtained from the normal direction to $Z_X$. As both $N_x Z_X$ and $N_y Z_Y$ are one-dimensional subspaces and $df_x$ gives a surjection, it is an isomorphism.
\end{rem}

We next discuss sections of $b$-maps, in preparation of Section \ref{sec:bhyperfibration}.
\begin{prop}\label{prop:bmapsections} Let $f:\, (X,Z_X) \to (Y,Z_Y)$ be a $b$-map, and let $s$ be a section of $f$. Then $s:\, (Y,Z_Y) \to (X,Z_X)$ is a $b$-map, and $\ker \b{df}_x \oplus \b{ds}_y(\b{T_y Y}) = \b{T_x X}$ for all $x \in X$, where $y = f(x)$.
\end{prop}
\bp By definition $f \circ s = {\rm id}_Y$, hence $s^{-1}(Z_X) = s^{-1}(f^{-1}(Z_Y)) = (f \circ s)^{-1}(Z_Y) = Z_Y$. Note that given finite-dimensional vector spaces $U, W, V', V$ with $V' \subset V$ so that $U \hookrightarrow V' \twoheadrightarrow W$ and $U \hookrightarrow V \twoheadrightarrow W$ we have $V' = V$ by counting dimensions. Let $y \in Z_Y$ and $x \in f^{-1}(y)$ be given. By definition $\ker df_x \hookrightarrow T_x X \twoheadrightarrow T_y Y$ using $df_x$. There is a clear surjection of $ds_y(T_y Y)$ onto $T_y Y$. As $ds_y(T_y Y) + T_x Z_X \subset T_x X$ and $\ker df_x \hookrightarrow T_x Z_X$ because $f^{-1}(Z_Y) = Z_X$, we see that $\ker df_x \hookrightarrow ds_y(T_y Y) + T_x Z_X \twoheadrightarrow T_y Y$, so that $ds_y(T_y Y) + T_x Z_X = T_x X$. We conclude that $s$ is a $b$-map. Similarly, consider $x \in X$ and denote $y = f(x)$. As $f$ is a $b$-map we have $\ker \b{df}_x \hookrightarrow \b{T_x X} \twoheadrightarrow \b{T_y Y}$ using $\b{df}_x$. Again there is a surjection of $\b{ds}_y(\b{T_y Y})$ onto $\b{T_y Y}$. Note that $\ker \b{df}_x \oplus \b{ds}_y(\b{T_y Y}) \subset \b{T_x X}$ and furthermore $\ker \b{df}_x \hookrightarrow \ker \b{df}_x \oplus \b{ds}_y \b{T_y Y} \twoheadrightarrow \b{T_y Y}$, so that $\ker \b{df}_x \oplus \b{ds}_y \b{T_y Y} = \b{T_x X}$.
\ep
Let $f:\, (X^4,Z_X) \to (\Sigma^2,Z_\Sigma)$ be a $b$-map for which $\ker \b{df}$ is even-dimensional. For example, this situation arises when there is a $b$-symplectic form $\btwo{\omega}_\Sigma$ on $(\Sigma,Z_\Sigma)$ and a \bacs $\btwo{J}$ on $(X,Z_X)$ so that $\btwo{J}$ is $(\btwo{\omega}_\Sigma, f)$-tame. Because $f$ is transverse to $Z_\Sigma$, we conclude that $Z_\Sigma$ cannot contain critical values of $f$. Namely, if $x$ were a critical point of $f$, by a dimension count and \autoref{prop:kerbdf} we would have $\ker df_x = T_x X$, hence $df_x(T_x X) + T_{f(x)} Z_Y = T_{f(x)} Z_Y$, not $T_{f(x)} Y$. Said differently, such a map $f:\, X^4 \to \Sigma^2$ can only be turned into a $b$-map if $Z_\Sigma$ is disjoint from the set of critical values of $f$.
\subsection{\Blog{s} on surfaces}
\label{sec:blogsurfaces}
In this section we discuss \blog{s} on a compact surface $\Sigma$. These were first studied in \cite{Radko02}, but we are interested in existence, i.e.\ which pairs $(\Sigma, Z_\Sigma)$ admit \blog{}s. Any surface $\Sigma$, compact or not, orientable or not, admits a \blog. This is because any vector bundle admits transverse sections, hence so does $\bigwedge^2 T\Sigma$; the Poisson condition is immediate in dimension two. However, not every pair $(\Sigma,Z_\Sigma)$ admits a \blog, with an easy counterexample being $(\R P^2, \emptyset)$. The hypersurface $Z_\Sigma$ must be the zero set $N(s)$ of some section $s \in \Gamma(\bigwedge^{2} T\Sigma)$, so that $Z_\Sigma$ represents $w_1(\bigwedge^2 T\Sigma)$ and $\Sigma \setminus Z_\Sigma$ is oriented by $s$. More precisely, we have the following.
\begin{prop}\label{prop:blogsurface} Let $(\Sigma^2, Z_\Sigma)$ be given. Then $(\Sigma, Z_\Sigma)$ admits a \blog if and only if $w_1(\Sigma) = {\rm PD}[Z_\Sigma]$, the Poincar\'e dual using $\Z_2$-coefficients.
\end{prop}
\bp The necessity of the condition follows after noting that $w_1(\bigwedge^2 T\Sigma) = w_1(T\Sigma) = w_1(\Sigma)$ and the \blog provides a transverse section of $\bigwedge^2 T\Sigma$ with zero set $Z_\Sigma$. To see sufficiency, then by hypothesis there exists a transverse section $\pi \in \Gamma(\bigwedge^2 T \Sigma)$ whose zero set $N(\pi)$ is equal to $Z_\Sigma$. As we are in dimension two, the Poisson condition $[\pi,\pi] = 0$ is immediate. Thus $\pi$ is a \blog with $Z_\pi = Z_\Sigma$.
\ep
If $\Sigma$ is orientable then $w_1(\Sigma) = 0$, making the condition that $Z_\Sigma$ must be $\Z_2$-nullhomologous. Remark further that $(\Sigma, Z_\Sigma)$ carries a \blog iff it is $b$-orientable.
%
%
\section{Construction of \blog{}s}
\label{sec:blogconstructions}
In this section we discuss extensions of symplectic results to the world of $b$-geometry to obtain existence results for \blog{}s. Let $X$ and $Y$ be compact connected manifolds, and assume that $Y$ is equipped with a \blog. Using \autoref{prop:blogbsymp} we can view $Y$ as being a $b$-manifold with a \bsymp, obtaining a triple $(Y,Z_Y,\btwo{\omega}_Y)$. Given a map $f:\, X \to Y$ such that $f$ is transverse to $Z_Y$, we can turn it into a $b$-map $f:\, (X,Z_X) \to (Y,Z_Y)$ by defining $Z_X := f^{-1}(Z_Y)$. We wish to equip $X$ with a \blog by constructing a $b$-symplectic form coming from $f$ and $\btwo{\omega}_Y$ and then use \autoref{prop:blogbsymp}. In analogy with \cite{Gompf04}, we will do this using \bacs{}s to guide the process.
\begin{defn}\label{defn:taming} Let $T:\, V \to W$ be a linear map between two finite-dimensional real vector spaces, and let $\omega$ be a skew-symmetric bilinear form on $W$. Then a linear complex structure $J$ on $V$ is called \emph{$(\omega, T)$-tame} if $T^* \omega(v, J v) > 0$ for all $v \in V \setminus \ker T$. For $f:\, (X,Z_X) \to (Y,Z_Y)$ a $b$-map and $\omega_Y \in \b{\Omega}^2(Y)$, a \bacs $J$ on $X$ is called \emph{$(\omega_Y, f)$-tame} if it is $(\omega_Y, \b{df}_x)$-tame for all $x \in X$.
\end{defn}
When $T = {\rm id}_V$ or $f = {\rm id}_X$ and $Z_X = \emptyset$, we recover the usual notion of $J$ being $\omega$-tame. Note that any $\omega \in \b{\Omega}^2(X)$ taming a \bacs $J$ is necessarily non-degenerate. Hence if $\omega \in \b{\Omega}^2(X)$ is a closed $b$-two-form taming some $J$, then $\omega$ is $b$-symplectic and $J$ induces the same $b$-orientation as $\omega$. Note moreover that if $J$ is $(\omega_Y, f)$-tame, then $\ker \b{df}$ is a $J$-complex subspace so that preimages of regular values of $f$ are $J$-holomorphic submanifolds. Indeed, if $v \in \ker \b{df}$ and $J v \not\in \ker \b{df}$, we would have $0 = f^*\omega(v, Jv) = f^*\omega(Jv, J(Jv)) > 0$, which is a contradiction. If $\omega$ is $b$-symplectic there exist tamed $J$ for $\omega$, and the space of such $J$ is convex and hence contractible \cite{McDuffSalamon98}.
\begin{rem}\label{rem:tamingopen} An important property of tameness is that the taming condition is open, i.e.\ it is preserved under sufficiently small perturbations of $\omega$ and $J$, and of varying the point in $X$. Namely, the taming condition $\omega(v, J v) > 0$ for the pair $(\omega, J)$ holds provided it holds for all $v \in \b{\Sigma X} \subset \b{TX}$, the unit sphere bundle with respect to some preassigned metric, as $X$ is compact. As $\b{\Sigma X}$ is also compact, the continuous map $\wt{\omega}:\, \b{\Sigma X} \to \R$ given by $\wt{\omega}(v) := \omega(v, J v)$ for $v \in \b{\Sigma X}$ is bounded from below by a positive constant on $\b{\Sigma X}$. But then this map will remain positive under small perturbations of $\omega$ and $J$. Similarly the condition of $\omega$ taming $J$ on $\ker \b{df}$ is open. Consider $x \in X$ so that $\omega(v, J v) > 0$ for all $v \in \ker \b{df}_x$. As $\wt{\omega}$ is continuous and $\b{\Sigma X}$ is compact, there exists a neighbourhood $U$ of $\ker \b{df}_x \cap \b{\Sigma X}$ in $\b{\Sigma X}$ on which $\wt{\omega}$ is positive. Points $x' \in X$ close to $x$ will then have $\ker \b{df}_{x'} \subset U$ because $\ker \b{df}$ is closed.
\end{rem}
We will not use the associated notion of compatibility, where $J$ also leaves $\omega$ invariant, as we use \bacs{}s as auxiliary structures to show non-degeneracy, and make use of the openness of this condition. For this reason, all \bacs{}s will only be required to be continuous as this avoids arguments to ensure smoothness. Using a \bacs on $X$, one can demand the existence of suitable forms on both the base and the fibers, to be combined into a $b$-symplectic form on the total space $X$. Note that the fibers here do not have to be two-dimensional.
\begin{rem} Products of log-symplectic manifolds are not log-symplectic in general, as the product form is not log-symplectic when both are bona fide. However, the product of a log-symplectic manifold with a symplectic manifold is always log-symplectic.
\end{rem}
Given a map $f:\, X \to Y$ between compact connected manifolds, we will use the following notation.
\bi
	 \item $F_y = f^{-1}(y)$ for $y \in Y$ is the level set, or \emph{fiber}, of $f$ over $y$;
	 \item $[F]$ is the homology class of a generic fiber, i.e.\ the inverse image of a regular value.
\ei

The latter will only be used in instances where it is well defined and independent of the regular value. We now formulate the $b$-version of \cite[Theorem 3.1]{Gompf04}. Note that there are no assumptions on the dimension of the manifolds involved, other than that they must be even. 
\begin{thm}\label{thm:bthurstontrick} Let $f:\, (X,Z_X) \to (Y,Z_Y)$ be a $b$-map between compact connected $b$-oriented $b$-manifolds, $\btwo{J}$ a \bacs on $(X,Z_X)$, $\btwo{\omega}_Y$ a $b$-symplectic form on $(Y,Z_Y)$, and $c \in H^2(X;\R)$. Assume that
\bi
	\item[i)] $\btwo{J}$ is $(\btwo{\omega}_Y,f)$-tame;
	\item[ii)] for each $y \in Y$, $F_y$ has a neighbourhood $W_y$ with a closed two-form $\eta_y \in \Omega^2(W_y)$ such that $[\eta_y] = c|_{W_y} \in H^2(W_y;\R)$, and $\rho_X^* \eta_y$ tames $\btwo{J}|_{\ker \b{df}_x}$ for all $x \in W_y$.
\ei
Then $(X,Z_X)$ admits a \blog.
\end{thm}
\begin{rem} As one can see by inspecting the proof, the second condition in the above theorem can be replaced by the existence of a closed two-form $\eta \in \Omega^2(X)$ such that $\rho_X^* \eta$ tames $\left.J\right|_{\ker \b{df}_x}$ for all $x \in X$. Such forms are called \emph{symplectic on the fibers of $f$}. Comparing \autoref{thm:bthurstontrick} to \cite[Theorem 3.1]{Gompf04} one sees that the form on the base $Y$ is $b$-symplectic, while the fiberwise form $\eta$ is still symplectic in the usual sense.
\end{rem}
The proof of this result is modelled on that by Gompf of \cite[Theorem 3.1]{Gompf04}.
\bp Using \autoref{prop:blogbsymp} it suffices to construct a $b$-symplectic form on the $b$-manifold $(X,Z_X)$. Let $\xi \in \Omega^2_{\rm cl}(X)$ be such that $[\xi] = c$. Then for each $y \in Y$ we have $[\eta_y] = c|_{W_y} = [\xi]|_{W_y}$, so on $W_y$ we have $\eta_y = \xi + d\alpha_y$ for some $\alpha_y \in \Omega^1(W_y)$. As each $X \setminus W_y$ and hence $f(X \setminus W_y)$ is compact, each $y \in Y$ has a neighbourhood disjoint from $f(X\setminus W_y)$. Cover $Y$ by a finite amount of such open sets $U_i$ so that each $f^{-1}(U_i)$ is contained in some $W_{y_i}$. Let $\{\varphi_i\}$ be a partition of unity subordinate to the cover $\{U_i\}$ of $Y$, so that $\{\varphi_i \circ f\}$ is  a partition of unity of $X$. Define a two-form $\eta \in \Omega^2(X)$ on $X$ via
\be
	\eta := \xi + d(\sum_i (\varphi_i \circ f) \alpha_{y_i}) = \xi + \sum_i (\varphi_i \circ f) d\alpha_{y_i} + \sum_i (d\varphi_i \circ df) \wedge \alpha_{y_i}.
\ee
Then $d\eta = 0$ and $[\eta] = c$. The last of the above three terms vanishes when applied to a pair of vectors in $\ker df_x$ for any $x \in X$, so on each $\ker df_x$ we have
\be
	\eta = \xi + \sum_i (\varphi_i \circ f) d\alpha_{y_i} = \sum_i (\varphi_i \circ f)(\xi + d\alpha_{y_i}) = \sum_i (\varphi_i \circ f) \eta_{y_i}.
\ee
Consider $\b{\eta} := \rho_X^* \eta$, a closed $b$-two-form on $(X,Z_X)$. By \autoref{prop:kerbdf} and the above we see that on $\ker \b{df}$, the $b$-form $\b{\eta}$ is a convex combination of $\btwo{J}$-taming $b$-forms, so $\btwo{J}|_{\ker \b{df}}$ is $\b{\eta}$-tame. Define for $t > 0$ a global closed $b$-two-form on $X$ by $\btwo{\omega}_t := f^* \btwo{\omega}_Y + t \b{\eta}$. We show that $\btwo{\omega}_t$ is $b$-symplectic by showing $\btwo{\omega}_t$ tames $\btwo{J}$ for $t$ small enough. By \autoref{rem:tamingopen} it is enough to show that there exists a $t_0 > 0$ so that $\btwo{\omega}_t(v, \btwo{J} v) > 0$ for every $t \in (0,t_0)$ and $v$ in the unit sphere bundle $\b{\Sigma X} \subset \b{TX}$ with respect to some metric. Note that for $v \in \b{TX}$ we have $\btwo{\omega}_t(v, \btwo{J} v) = f^* \btwo{\omega}_Y(v, \btwo{J} v) + t \b{\eta}(v, \btwo{J} v)$. As $\btwo{J}$ is $(\btwo{\omega}_Y,f)$-tame, the first term is positive for $v \in \b{TX} \setminus \ker \b{df}$ and is zero otherwise. The second term $\b{\eta}(v, \btwo{J} v)$ is positive on $\ker \b{df}$ because $\btwo{J}|_{\ker \b{df}}$ is $\b{\eta}$-tame, hence is also positive for all $v$ in some neighbourhood $U$ of $\ker \b{df} \cap \b{\Sigma X}$ in $\b{\Sigma X}$ by openness of the taming condition. We conclude that $\btwo{\omega}_t(v, \btwo{J} v) > 0$ for all $t > 0$ when $v \in U$. The function $\b{\eta}(v, \btwo{J} v)$ is bounded on the compact set $\b{\Sigma X} \setminus U$. Furthermore, $f^*\btwo{\omega}_Y(v, \btwo{J} v)$ is bounded from below there by a positive constant, as it is positive away from $\ker \b{df}$, and thus also away from $\ker \b{df} \cap \b{\Sigma X} \subset U$. But then $\btwo{\omega}_t(v, \btwo{J} v) > 0$ for all $0 < t < t_0$ for $t_0$ sufficiently small, so that $\btwo{\omega}_t$ is indeed $b$-symplectic for $t$ small enough.
\ep
\begin{rem} Given a $b$-map $f:\, (X,Z_X) \to (Y,Z_Y)$, the singular locus of $X$ is given by $Z_X = f^{-1}(Z_Y)$, so that it consists purely of fibers. This means the fibers of $f$, which are natural candidates for being $b$-symplectic submanifolds of $X$, never hit $Z_X$.
\end{rem}
With \autoref{thm:bthurstontrick} in hand we can prove the $b$-version of Thurston's result for symplectic fiber bundles with two-dimensional fibers \cite{Thurston76}, adapting the proof by Gompf in \cite{Gompf01}.
\begin{thm}\label{thm:bthurstonfibration} Let $f:\, (X^{2n},Z_X) \to (Y^{2n-2}, Z_Y)$ be a $b$-fibration between compact connected $b$-manifolds. Assume that $(Y,Z_Y)$ is $b$-symplectic and that the generic fiber $F$ is orientable and $[F] \neq 0 \in H_2(X;\R)$. Then $(X,Z_X)$ admits a \blog.
\end{thm}
This theorem has the following immediate corollary, phrased without using the $b$-language.
\begin{cor}\label{cor:fibrationblogs} Let $f:\, X^{2n} \to Y^{2n-2}$ be a fibration between compact connected manifolds. Assume that $Y$ admits a \blog $\pi$ and that the generic fiber $F$ is orientable and $[F] \neq 0 \in H_2(X;\R)$. Then $(X,Z_X)$ admits a \blog for $Z_X = f^{-1}(Z_{\pi})$.
\end{cor}
\bp Let $Z_Y := Z_{\pi}$. Then by \autoref{prop:blogbsymp}, $(Y,Z_Y)$ is $b$-symplectic. By \autoref{rem:fibrintobfibr}, the map $f$ can be seen as a $b$-fibration $f:\, (X,Z_X) \to (Y,Z_Y)$, where $Z_X = f^{-1}(Z_Y)$. Using \autoref{thm:bthurstonfibration} we conclude that $(X,Z_X)$ admits a \blog.
\ep
For example, any homologically essential oriented surface bundle over a surface is log-symplectic (regardless of whether the base is orientable). The fact that any such base surface is log-symplectic is discussed in Section \ref{sec:blogsurfaces}.
\bp[ of \autoref{thm:bthurstonfibration}] We check the conditions of \autoref{thm:bthurstontrick}. If necessary, pass to a finite cover of $Y$ so that the fibers of $f$ are connected. Let $\btwo{\omega}_Y$ be a $b$-symplectic form on $Y$ and choose an $\omega_Y$-taming \bacs $\btwo{J}_Y$ on $Y$. Fix an orientation for the generic fiber $F$, which orients $\ker df$ and $\ker \b{df}$, using \autoref{prop:kerbdf}. As $[F] \neq 0$, this forces $f:\, (X,Z_X) \to (Y,Z_Y)$ to be an orientable $b$-fibration. In particular, $f$ provides a $b$-orientation for $(X,Z_X)$ by using the fiber-first convention. Let $g$ be a metric on $\b{TX}$ and let $H \subset \b{TX}$ be the subbundle of orthogonal complements to $\ker \b{df}$. Note here that $\ker \b{df}$ is the tangent space to the fibers of $f$, because $f$ is a $b$-fibration. Then $\b{df}|_H:\, H \to \b{TY}$ is an isomorphism on each fiber. Define a \bacs $\btwo{J}$ on $X$ as follows. On $H$, define $\btwo{J}|_H = f^* \btwo{J}_Y$ using the isomorphism given by $\b{df}$. On the two-dimensional tangent spaces to the fibers, use the metric and define $\btwo{J}$ by $\frac{\pi}{2}$-counterclockwise rotation, demanding that $f$ is orientation preserving via the fiber-first convention. This determines $\btwo{J}$ uniquely on $\b{TX}$ by linearity. Moreover, $\btwo{J}$ is $(\btwo{\omega}_Y, f)$-tame as $f^* \btwo{\omega}_Y(v, \btwo{J} v) = \btwo{\omega}_Y(\b{df} v, \btwo{J}_Y \b{df} v) > 0$, for all $v \in \b{TX} \setminus \ker \b{df} \cong H$.

We construct the required neighbourhoods $W_y$ and forms $\eta_y$ for each fiber. Let $c \in H^2(X;\R)$ be through duality a class such that $\langle [F], c\rangle = 1$, using that $[F] \neq 0$. Given $y \in Y$, let $D_y \subset Y$ be an open disk containing $y$, fully contained in a trivializing neighbourhood of $f$ around $y$. Define $W_y := f^{-1}(D_y) \cong D_y \times F_y$. Using that $F_y$ is two-dimensional, choose an area form on $F_y$ inducing the preimage orientation of the fiber, and let $\eta_y \in \Omega^2(W_y)$ be the pullback of this form via the projection $p:\, W_y \to F_y$. Because $\langle [F_y], \eta_y \rangle = 1 = \langle [F_y], c\rangle$ and $H_2(W_y;\R)$ is generated by $[F_y]$,  it follows that $[\eta_y] = c|_{W_y} \in H^2(W_y;\R)$. To check that $\rho_X^* \eta_y$ tames $\btwo{J}$ on $\ker \b{df}_x$ for $x \in W_y$, recall that there $\btwo{J}$ is defined via rotation. As $\rho_X^* \eta_y$ is the pullback of the area form of a fiber, taming follows as its restriction to a fiber is an area form for that fiber. Applying \autoref{thm:bthurstontrick} we obtain a \blog on $(X,Z_X)$.
\ep
\begin{rem} Note how the two-dimensionality of the fibers is used to be able to define an \acs there, and for choosing an area form.
\end{rem}
%
%
\section{\Alf{}s}
\label{sec:achirallefschetz}
In this section we revisit the result proven in \cite{Cavalcanti13} that any homologically essential \alf with orientable fibers gives rise to a \blog in dimension four, which can be chosen to be bona fide. Our proof using $b$-geometry is a direct adaptation of the proof that four-dimensional \lf{}s give rise to symplectic structures, phrased using \acs{}s.
\begin{defn}\label{defn:achirallefschetzfibration} An \emph{\alf} is a map $f:\, X^{2n} \to \Sigma^2$ between compact connected manifolds so that for each critical point $x \in X$ there exist complex coordinate charts centered at $x$ and $f(x)$ in which $f$ takes the form $f:\, \C^n \to \C$, $f(z_1, \dots, z_n) = z_1^2 + \dots + z_n^2$.
\end{defn}
Note that we do not require $X$ nor $\Sigma$ to be orientable. If they are however, after choosing orientations one can assign a sign to each critical point of $f$. A given critical point $x \in X$ obtains a sign by demanding that the complex structure of the chart on $\Sigma$ is compatible with its orientation; we then say $x$ is \emph{positive} if the complex structure on $X$ is compatible with its orientation, and \emph{negative} otherwise. Note that any \lf is also an \alf.
\begin{thm}[{\cite[Theorem 6.7]{Cavalcanti13}}] \label{thm:alflogsymp} Let $f:\, X^4 \to \Sigma^2$ be an \alf between compact connected manifolds. Assume that the generic fiber $F$ is orientable and $[F] \neq 0 \in H_2(X;\R)$. Then $X$ admits a \blog.
\end{thm}
This theorem should be viewed as a direct analogue of \cite[Theorem 10.2.18]{GompfStipsicz99} that homologically essential four-dimensional Lefschetz fibrations provide symplectic structures. We prove it by first showing that an \alf gives rise to what we call a \emph{\blf}. Noting \autoref{prop:kerbdf} we use the following notation, given a map $f:\, X \to Y$ between manifolds (in what follows, $Y = \Sigma$).
\bi
	 \item $\Delta = {\rm Crit}(f) \subset X$ is the set of critical points of $f$;
	 \item $\Delta_y = \Delta \cap F_y$ for $y \in Y$ is the set of critical points of $f$ lying on the fiber $F_y$;
	 \item $\Delta' = {\rm Sing}(f) \subset Y$ is the set of singular values of $f$.
\ei
\begin{defn}\label{defn:blefschetzfibration} A \emph{\blf} is a $b$-map $f:\, (X^{2n},Z_X) \to (\Sigma^2, Z_\Sigma)$ between compact connected $b$-oriented $b$-manifolds so that for each critical point $x \in \Delta$ there exist complex coordinate charts compatible with orientations induced from the $b$-orientations centered at $x$ and $f(x)$ in which $f$ takes the form $f:\, \C^n \to \C$, $f(z_1, \dots, z_n) = z_1^2 + \dots + z_n^2$.
\end{defn}
\begin{rem}
Given a \blf $f:\, (X,Z_X) \to (\Sigma,Z_\Sigma)$, the local model for $f$ around critical points $x \in \Delta$ implies $\ker df_x = T_x X$. Because $f$ is a $b$-map so is transverse to $Z_\Sigma$, we conclude that $Z_\Sigma$ and $\Delta'$ are disjoint. The $b$-orientation induces an orientation away from the singular locus so that it makes sense to demand compatibility of the charts.
\end{rem}

Alternatively, we could first define the notion of a $b$-\alf and then note that its critical values must be disjoint from the singular locus, so that we can then further demand compatibility of the charts specifying the local model of $f$.
\begin{prop}\label{prop:alfblf} Let $f:\, X^{2n} \to \Sigma^2$ be an \alf between compact connected manifolds which is injective on critical points. Assume that the generic fiber $F$ is orientable and $[F] \neq 0 \in H_{2n-2}(X;\R)$. Then there exists a hypersurface $Z_\Sigma \subset \Sigma$ so that $f:\, (X^{2n},Z_X) \to (\Sigma^2,Z_\Sigma)$ is a \blf, where $Z_X = f^{-1}(Z_\Sigma)$.
\end{prop}
\bp We deal with orientations as in \cite[Theorem 6.7]{Cavalcanti13}. Fix an orientation for the generic fiber $F$. As $[F] \neq 0$, this forces $f:\, X \setminus \Delta \to \Sigma \setminus \Delta'$ to be an orientable fibration, which in turn orients all fibers, including the singular ones. We conclude that $X$ is orientable if and only if $\Sigma$ is. If they are, choose orientations and split $\Delta'$ into disjoint sets $\Delta'_+$ and $\Delta'_-$ according to the sign of the critical points. Then, pick a separating curve $\gamma \subset \Sigma$ disjoint from $\Delta'$ such that its interior contains all of $\Delta'_-$ and no points from $\Delta'_+$. If $\Sigma$ is not orientable, there instead exists a curve $\gamma \subset \Sigma$ so that $\Sigma \setminus \gamma$ is orientable, hence so is $X \setminus f^{-1}(\gamma)$. Choose orientations and then homotope $\gamma$ through negative critical values of $f$ so that all critical points are positive. Define $Z_\Sigma := \gamma$ and let $Z_X := f^{-1}(Z_\Sigma)$. Because $Z_\Sigma$ does not hit $\Delta'$, it is immediate that $f$ is a $b$-map from $(X,Z_X)$ to $(Y,Z_\Sigma)$. Moreover, the orientations we chose give the appropriate $b$-orientations. But then $f$ is a \blf.
\ep
\begin{rem}\label{rem:curvechoice} The curve $\gamma$ used in the previous proof is not unique. For example, in the orientable case we chose a separating curve, but we could just as well have chosen disjoint curves around each negative critical point separately. The effect of this is that the $b$-manifold structures that are used are not unique either.
\end{rem}
Given a four-dimensional \blf, the $b$-version of Gompf's proof of \cite[Theorem 10.2.18]{GompfStipsicz99} in terms of tamed \acs{}s supplies us with a $b$-symplectic form on $(X,Z_X)$, giving a \blog on $X$.
\begin{thm}\label{thm:blfblog} Let $f:\, (X^4,Z_X) \to (\Sigma^2,Z_\Sigma)$ be a \blf between compact connected $b$-oriented $b$-manifolds which is injective on critical points. Assume that the generic fiber $F$ is orientable and $[F] \neq 0 \in H_2(X;\R)$. Then $(X,Z_X)$ admits a \blog.
\end{thm}
\bp Note that the exact sequence in homotopy, $\pi_1(F) \to \pi_1(X) \to \pi_1(Y) \to \pi_0(F) \to 0$, which was used in the proof of \autoref{thm:bthurstonfibration}, is also available for \lf{}s by \cite[Proposition 8.1.9]{GompfStipsicz99} because neighbourhoods of critical values admit specific sections. This still holds for \alf{}s and hence also for \blf{}s.  If necessary, use this exact sequence in homotopy to lift $f$ to a cover of $\Sigma$ so that $f$ has connected fibers. Fix an orientation for the generic fiber $F$, which orients $\ker df$ and $\ker \b{df}$ at regular points, using \autoref{prop:kerbdf}. As $[F] \neq 0$, this forces $f:\, (X \setminus f^{-1}(\Delta'),Z_X \setminus f^{-1}(\Delta')) \to (Y\setminus \Delta',Z_Y\setminus \Delta')$ to be an orientable $b$-fibration. Let $c \in H^2(X;\R)$ be a class dual to $[F] \in H_2(X;\R)$, i.e.\ such that $\langle c, [F] \rangle = 1$. Let $\b{\pi} \in \Gamma(\b{T\Sigma})$ be a transverse section specifying the $b$-orientation of $(\Sigma,Z_\Sigma)$. Then by \autoref{prop:blogbsymp}, $\rho_\Sigma(\b{\pi})$ is a \blog on $(X,Z_\Sigma)$, and $\btwo{\omega}_\Sigma := (\b{\pi})^{-1}$ is a $b$-symplectic form on $(\Sigma,Z_\Sigma)$.

We show the conditions of \autoref{thm:bthurstontrick} are satisfied. Use the proof of \autoref{thm:bthurstonfibration} to obtain a $(\btwo{\omega}_\Sigma, f)$-tame \bacs $\btwo{J}$ on $\b{TX}|_{X \setminus f^{-1}(\Delta')}$. As $\Delta'$ is disjoint from $Z_\Sigma$, let $V \subset \Sigma$ be the disjoint union of open balls $V_y$ disjoint from $Z_\Sigma$ and centered at each point $y \in \Delta'$. Set $W := f^{-1}(V) \subset X$ to be the union of the neighbourhoods $W_y := f^{-1}(V_y)$ of singular fibers $F_y$. Let $C \subset X$ be the disjoint union of open balls $C_y$ centered at each point $f^{-1}(y) = x \in \Delta_y$ for all $y \in \Delta'$ disjoint from $Z_X$ so that on each ball $f$ takes on the local form in \autoref{defn:blefschetzfibration}. Possibly shrink $C$ so that $\overline{C_y} \subset W_y$. The local description of $f$ gives an \acs on $C$ with the fibers being holomorphic, and we glue this to the existing $\btwo{J}$ on $X \setminus C$ using a trivial version of \autoref{prop:acsconstruct} to obtain a global $(\btwo{\omega}_\Sigma, f)$-tame \bacs $\btwo{J}$ on $(X,Z_X)$.

Let $y \in Y$ be given. If $y \in Y \setminus \Delta'$, proceed as in the proof of \autoref{thm:bthurstonfibration} to obtain a neighbourhood $W_y$ of $F_y$ and a closed two-form $\eta_y \in \Omega^2(W_y)$ such that $[\eta_y] = c|_{W_y}$ and such that $\eta_y$ tames $\btwo{J}$ on $\ker \b{df}_x$ for all $x \in W_y$. If $y \in \Delta'$, the singular fiber $F_y$ either is indecomposable or consists of two irreducible components $F_y^\pm$ which satisfy $[F_y^+] \cdot [F_y^-] = 1$ and $[F_y^\pm]^2 = -1$, see \cite{GompfStipsicz99}. In the latter case, note that $0 < 1 = \langle c, [F] \rangle = \langle c, [F_y]\rangle = \langle c, [F_y^+]\rangle + \langle c, [F_y^-]\rangle$. If either term is nonpositive assume without loss of generality that $\langle c, [F_y^-]\rangle = r \leq 0$. Define $c' := c + (\frac12 - r) c_y^+$, where $c_y^+ \in H^2(X;\R)$ is a class dual to $[F_y^+]$. As $[F_y] \cdot [F_y^\pm] = 0$ we then have $\langle c', [F] \rangle = \langle c, [F] \rangle > 0$, and furthermore $\langle c', [F_y^+]\rangle = \langle c, [F_y^+] \rangle +(\frac12 - r) > 0$ and $\langle c', [F_y^-]\rangle = \frac12 > 0$. Moreover, as different fibers do not intersect, we have $c|_{W_{y'}} = c'|_{W_{y'}}$ for $y' \neq y$. After finitely many repetitions, at most once for each $y \in \Delta'$, one obtains a class, again denoted by $c$, pairing positively with every fiber component (see \cite[Exercise 10.2.19]{GompfStipsicz99}).

Return to $y \in \Delta'$ and let $\sigma$ be the standard symplectic form on $C_y$ given locally in real coordinates by $\sigma = dx_1 \wedge dy_1 + dx_2 \wedge dy_2$, where $z_i = x_i + i y_i$. As all fibers $F_y'$ in $C_y$ are holomorphic, $\rho_X^* \sigma|_{F_{y'} \cap C_y}$ tames ${}^b J$ for all $y' \in f(C_y)$, so that $\rho_X^* \sigma$ tames $\btwo{J}$ on $C_y$. Let $\sigma_y$ be an extension of $\sigma$ to $F_y$ as a positive area form with total area $\langle \sigma_y, [F_y] \rangle$ equal to $\langle c, [F_y] \rangle$. Let $p:\, W_y \to F_y$ be a retraction and let $f:\, C_y \to [0,1]$ be a smooth radial function so that $f \equiv 0$ in a neighbourhood of $x$ and $f \equiv 1$ in a neighbourhood of $\partial C_y$, which is smoothly extended to $W_y$ by being identically $1$ outside $C_y$. On the ball $C_y$, the form $\sigma$ is exact, say equal to $\sigma = d\alpha$ for $\alpha \in \Omega^1(C_y)$. Define a two-form $\eta_y$ on $W_y$ by $\eta_y := p^*(f\sigma_y) + d((1-f)\alpha)$, which is closed as $f\sigma_y$ is a closed area form on $F_y$. Near $x$ we have $f \equiv 0$ so that $\eta_y = d\alpha = \sigma$, hence there $\rho_X^* \eta_y = \rho_X^* \sigma$ tames $\btwo{J}$, hence in particular tames $\btwo{J}|_{\ker \b{df}}$. Similarly, $\sigma_y$ is an area form on $F_y \setminus \{x\}$ for the orientation given by $\btwo{J}$, using \autoref{prop:kerbdf}. But then $\rho_X^* \sigma_y$ tames $\btwo{J}$ on $TF_y = \ker df \cong \ker \b{df}$ on $F_y$, so that the same holds for $\rho_X^* \eta_y$ as this condition is convex. By openness of the taming condition, shrinking $V_y$ and hence $W_y$ and possibly $C_y$ we can ensure that $\rho_X^* \eta_y$ tames $\btwo{J}|_{\ker \b{df}}$ on $W_y$. Finally, note that $[\eta_y] = c|_{W_y} \in H^2(W_y;\R)$ by construction. We have now obtained the required neighbourhoods $W_y$ and forms $\eta_y$ for all $y \in Y$ so that by \autoref{thm:bthurstontrick} we obtain a \blog on $(X,Z_X)$.
\ep
\bp[ of \autoref{thm:alflogsymp}] This follows as an immediate corollary to \autoref{prop:alfblf} and \autoref{thm:blfblog}. We require that $f$ is injective on critical points so that critical values also obtain a sign, allowing us to group them appropriately. Here we note that we can always perturb $f$ so that it is injective on critical points. If one does not want to assume this one can proceed as in \cite[Lemma 3.3]{Gompf05}.
\ep
\begin{rem} The \blog constructed in the proof of \autoref{thm:alflogsymp} has connected singular locus, and may be empty. By \autoref{rem:curvechoice} we could have easily ensured that the singular locus has multiple components. Moreover, by the following theorem we can add components to the singular locus at will using a neighbourhood of a fiber. In particular we can always ensure the \blog is bona fide.
\end{rem}
\begin{thm}[{\cite[Theorem 5.1]{Cavalcanti13}}] \label{thm:addsinglocus} Let $(X^{2n},Z_X,\pi)$ be a compact log-symplectic manifold and let $k > 0$ be an integer. Assume that $X$ has a compact symplectic submanifold $F^{2n-2} \subset X \setminus Z_X$ with trivial normal bundle. Then $(X,Z_X')$ admits a \blog $\pi'$ agreeing with $\pi$ away from $Z_X' \setminus Z_X$, where $Z_X'$ is the disjoint union of $Z_X$ with $k$ copies of $F \times S^1$.
\end{thm}
It is an interesting question whether every \blog on a four-manifold can be obtained out of an \alf using the construction of \autoref{thm:alflogsymp}. This parallels the development by Gompf and Donaldson between Lefschetz pencils and symplectic structures on four-manifolds. The first thing to note is that in our construction the fibers are always contained in the eventual singular locus, or are disjoint from it. Moreover, all \blog{}s we construct are \emph{proper}, in that all connected components of their singular loci are compact and have a compact symplectic leaf. This implies compact log-symplectic four-manifolds $(X,Z_X,\pi)$ for which $\pi$ is not proper are not reached by our construction. Note however that one can deform $\pi$ to a proper \blog if all components of $Z_X$ are compact \cite[Theorem 3.6]{Cavalcanti13}. More serious is the fact that in our construction the singular locus fibers over a circle in the base with specified diffeomorphism type of the fiber. Hence our construction cannot create log-symplectic four-manifolds $X$ with a disconnected singular locus $Z_X$ with at least two components not fibering over $S^1$ with the same genus fiber. The following example shows this can indeed happen.
\begin{exa} Let $X^4 = \Sigma_g \times \Sigma_h$ be the product of compact surfaces of genera $g \geq 2$ and $h \geq 1$ carrying the product symplectic form, and consider the map $f:\, X \to \Sigma_h$ given by projection. Consider a copy of the torus $T^2 = S^1 \times S^1$ by picking an essential circle in both the base and the fiber. It is Lagrangian and homologically nontrivial, so that by a result of Gompf \cite[Lemma 1.6]{Gompf94}, by a slight perturbation of the symplectic structure in a neighbourhood of the torus we find a symplectic structure on $X$ for which $T^2$ is symplectic. Applying \autoref{thm:addsinglocus} we obtain a log-symplectic structure on $(X, T^3)$. Now use a fiber $\Sigma_g$ of $f$ disjoint from the torus and apply \autoref{thm:addsinglocus} once more to obtain a \blog on $(X, Z_X)$, where $Z_X = T^3 \sqcup \Sigma_g \times S^1$. 

Note that $T^3$ cannot fiber over $S^1$ with fibers of genus other than one. Any fibration $p:\, T^3 \to S^1$ with fiber $F$ induces a long exact sequence in homotopy groups, a part of which reads $\pi_2(S^1) \to \pi_1(F) \to \pi_1(T^3) \to \pi_1(S^1)$, or more concretely $0 \to \pi_1(F) \hookrightarrow \Z^3 \to \Z$. This shows that $\pi_1(F)$ injects into the Abelian group $\Z^3$, hence must itself be Abelian and cannot have torsion. By counting its rank must be two, so that that the genus of $F$ must be one, and the fibers of $p$ are tori. Similarly, the product $Y := \Sigma_g \times S^1$ cannot fiber over $S^1$ with fiber $F$ being a torus. We have $b_1(Y) = 1 + b_1(\Sigma_g) = 2g +1 \geq 5$ as $g \geq 2$. However, if $F$ were a torus, $Y$ would be a mapping torus of $T^2$, hence $b_1(Y) \leq 1 + b_1(T^2) = 3$ which is a contradiction. By the discussion preceding this example we conclude that the log-symplectic structure on $(X,Z_X)$ cannot be obtained from the construction in \autoref{thm:alflogsymp}.
\end{exa}
%
%
%
\section{\texorpdfstring{$b$}{b}-Hyperfibrations}
\label{sec:bhyperfibration}
In this section we introduce a class of fibration-like maps with two-dimensional fibers that can be made to satisfy the conditions of \autoref{thm:bthurstontrick}. We call them \emph{$b$-hyperfibrations}, as they are the $b$-analogue of the notion of a hyperpencil with empty base locus, as introduced by Gompf in \cite{Gompf04}. After defining them we show that a $b$-hyperfibration satisfying a condition analogous to having homologically essential fibers gives rise to a \blog.

Let $E, F \to X$ be real vector bundles over a manifold $X$ and let $T:\, E \to F$ be a continuous bundle map. Call a point $x \in X$ \emph{regular} if $T_x:\, E_x \to F_x$ is surjective, and \emph{critical} otherwise. Let ${\rm Reg}(T)$ denote the spaces of regular points of $T$. Form the space $P \subset E$ by
\be
	P = \ol{\bigcup\limits_{x \in {\rm Reg}(T)} \ker T_x},
\ee
and let $P_x = P \cap E_x$ for $x \in X$. Then $P_x = \ker T_x$ when $x$ is regular, and otherwise $P_x \subset \ker T_x$ consists of all limits of sequences of vectors at regular points which are annihilated by $T$.
\begin{defn}\label{defn:wrapped} A point $x \in X$ is called \emph{$T$-wrapped} if ${\rm span}_\R P_x$ has real codimension at most two in $\ker T_x$.
\end{defn}
Note that all points $x \in {\rm Reg}(T)$ are wrapped as then $P_x = \ker T_x$. In our applications $\ker T_x$ will be even-dimensional, so that the wrappedness condition is immediate unless ${\rm rank}(E) \geq 6$. Let $f:\, (X,Z_X) \to (Y,Z_Y)$ be a $b$-map between compact connected $b$-manifolds.  As before, let $\Delta$ and $\Delta_y$ be the set of critical points of the map $f:\, X \to Y$ and those on the fiber $F_y$ for $y \in Y$ respectively. We can then apply the above definition to the continuous or in fact smooth map $\b{d} f:\, \b{TX} \to \b{TY}$, noting \autoref{prop:kerbdf}. Note that if $X$ is four-dimensional and $Y$ is a surface, every critical point of $f$ will be automatically $\b{df}$-wrapped if $\ker \b{df}$ is always even-dimensional. This is because $\ker \b{df}$ is two-dimensional at regular points and the dimension of $\ker \b{df}$ cannot exceed four at singular points. With this we can give the definition of a $b$-hyperfibration.
\begin{defn}\label{defn:bhyperfibration} A \emph{$b$-hyperfibration} is a $b$-map $f:\, (X^{2n},Z_X) \to (Y^{2n-2},Z_Y)$ between compact connected $b$-oriented $b$-manifolds so that there exists a \blog $\btwo{\omega}_Y$ on $(Y, Z_Y)$ and such that
\bi
	\item[i)] each critical point $x \in \Delta$ is $\b{df}$-wrapped;
	\item[ii)] for each critical point $x \in \Delta$, there exists a neighbourhood of $x$ and a \bacs $\btwo{J}_x$ on the $b$-manifold $(W_x, W_x \cap Z_X)$ such that $\btwo{J}_x$ is $(\btwo{\omega}_Y,f)$-tame;
	\item[iii)] for each $y \in Y$ there are only finitely many critical points $\Delta_y$ lying on its fiber $F_y$.
\ei
\end{defn}
This should be compared with \cite[Definition 2.4]{Gompf04}. Note that the definition does not require $X$ nor $Y$ to be orientable. When $Z_Y = \emptyset$ one almost recovers the definition by Gompf of a hyperpencil with empty base locus, over an arbitrary symplectic base.
\begin{rem} While a $b$-hyperfibration $f$ may have infinitely many critical points, note that regular points of $f$ are dense in $X$, arguing as in \cite[Theorem 2.11]{Gompf04}. If an open $W \subset X$ would consist entirely of critical points, choose a point $x_0 \in W$ which minimizes $\dim \ker \b{df}_x$, and using \autoref{prop:kerbdf} note that $\ker \b{df} \cong \ker df$ is a smooth distribution near $x_0$ as it can be realized as $\ker d(\pi \circ f)$ for a projection $\pi$. Then take a vector field in $\ker df$ and integrate it to obtain a curve of critical points all lying in a single fiber of $f$. This contradicts assumption iii) in the definition of a $b$-hyperfibration.
\end{rem}
The wrappedness of the critical points will be used to obtain a global $(\btwo{\omega}_Y,f)$-tame \bacs out of the locally existing ones. With this new notion in hand we can move on to the following result, which is the appropriate $b$-analogue of \cite[Theorem 2.11]{Gompf04}. Given a $b$-hyperfibration $f:\, (X,Z_X) \to (Y,Z_Y)$ and $y \in Y$ fixed, we refer to the closures of the connected components of $F_y \setminus \Delta_y$ as the \emph{components} of the fiber $F_y$.
\begin{thm}\label{thm:bhyperfibration} Let $f:\, (X,Z_X) \to (Y, Z_Y, \btwo{\omega}_Y)$ be a $b$-hyperfibration between compact connected $b$-oriented $b$-manifolds. Assume that there exists a finite collection $S$ of sections of $f$ interescting all fiber components non-negatively and for each fiber component at least one section in $S$ intersecting positively. Then $(X,Z_X)$ admits a \blog.
\end{thm}
Note that the condition on the existence of such a collection $S$ of fibers of $f$ implies that each component of each fiber is homologically essential.
\begin{rem} In the above theorem, note that by \autoref{prop:bmapsections} sections $s \in S$ are $b$-maps $s:\, (Y,Z_Y) \to (X,Z_X)$ and furthermore that $\ker \b{df}_x \oplus \b{ds}_y(\b{T_y Y}) = \b{T_x X}$ for all $x \in X$, where $y = f(x)$. In the proof of \autoref{thm:bhyperfibration} we show a $b$-hyperfibration naturally gives rise to a \bacs which is $(\btwo{\omega}_Y,f)$-tame so that $\ker \b{df}$ carries a $b$-orientation. By \autoref{prop:kerbdf} for smooth points $x \in X \setminus \Delta$ we then have an orientation for $T_x F_y = \ker df_x \cong \ker \b{df}_x$. Both $\b{T_x X}$ and $\b{T_y Y}$ carry orientations, hence so does $\b{ds}_y(\b{T_y Y})$. We can then define the positive intersection of $s \in S$ with $F_y$ by comparing the $b$-orientations on these tangent spaces in the usual way. Note that $s$ must intersect fibers in smooth points of $f$ as it is a section.
\end{rem}
The proof of \autoref{thm:bhyperfibration} will be modelled on Gompf's proof of \cite[Theorem 2.11]{Gompf04}. There will be an interplay between two types of singular behavior, namely that of the $b$-manifold structure and that of the fibration itself. The relation between these has been discussed before in \autoref{prop:kerbdf}. Note that in the case of a $b$-hyperfibration $f:\, (X,Z_X) \to (Y, Z_Y)$ with $Z_Y = \emptyset$, we cannot apply \cite[Theorem 2.11]{Gompf04} directly. Indeed, there is no base locus, but instead a set of sections $S$. This is akin to obtaining a \lf out of a Lefschetz pencil by blowing up the base locus (see also \autoref{rem:homessential}).

For $V$ a real finite-dimensional vector space, let $\mc{B}_V \subset {\rm Aut}(V)$ be the open set of linear operators on $V$ with no real eigenvalues, and $\mc{J}_V \subset \mc{B}_V$ the set of complex structures on $V$ in either orientation. The following lemma is proven in \cite{Gompf04} as Corollary 4.2, loc.\ cit.
\begin{lem}\label{lem:acsretraction} Let $V, W$ be real finite-dimensional vector spaces. Then there exists a canonical real-analytic retraction $j:\, \mc{B}_V \to \mc{J}_V$, satisfying for all linear maps $T:\, V \to W$ such that $T A = B T$, that $T j(A) = j(B) T$ (whenever both sides are defined).
\end{lem}
Because the retraction in the previous lemma is canonical we can apply it pointwise to a continuously varying map, to again obtain a continuous map. Let $E^{2n}, F^{2n-2} \to X$ be real oriented vector bundles over a compact manifold $X$. In what follows, a \emph{two-form} on a vector bundle is a continuously varying choice of skew-symmetric bilinear form on each fiber. The next proposition can be extracted from \cite[Lemma 3.2]{Gompf04}. We include a proof for completeness.
\begin{prop}\label{prop:acsconstruct} Let $T:\, E \to F$ be a continuous bundle map and $\omega_F$ a non-degenerate two-form on $F$. Assume that for all $x \in X$ there exists a neighbourhood $W_x$ of $x$ with an $(\omega_F, T)$-tame complex structure on $E|_{W_x}$. Assume that each critical point $x \in {\rm Crit}(T)$ is wrapped. Then there exists a continuous $(\omega_F, T)$-tame complex structure $J$ on $E$.
\end{prop}
\bp
Cover $X$ by open sets $W_\alpha$ equipped with complex structures $J_\alpha$ on $E|_{W_\alpha}$ which are $(\omega_F, T)$-tame. Let $\{\varphi_\alpha\}$ be a subordinate partition of unity, and define
\be
	A := \sum_\alpha \varphi_\alpha J_\alpha:\, E \to E, \qquad\qquad B := \sum_\alpha \varphi_\alpha T_* J_\alpha:\, T(E) \to T(E),
\ee
so that $TA = BT$. Since $\ker T$ is $J_\alpha$-complex for all $\alpha$ because $J_\alpha$ is $(\omega_F, T)$-tame (see below \autoref{defn:taming}), $J_\alpha$ descends to a map $T_* J_\alpha:\, T(E) \to T(E)$, hence $B$ is well-defined. In order to apply \autoref{lem:acsretraction} we show that $A_x \in \mc{B}_{E_x}$ for all $x \in X$, i.e.\ that $A$ has no real eigenvalues. Let $\lambda$ be an eigenvalue of $A$ with eigenvector $v \in E$. Then $B T v = T A v = T \lambda v = \lambda T v$, so either $T v = 0$, or $T v$ is a $\lambda$-eigenvector for $B$. As each $T_* J_\alpha$ is $\omega_F$-tame, $B$ has no real eigenvalues on any fiber. Indeed, for $0 \neq w \in T(E)$, we have $\omega_F(w,w) = 0$. Hence if $B w = \lambda w$ for some $\lambda \in \R$, we have
\be
	0 = \omega_F(w, \lambda w) = \omega_F(w, B w) = \omega_F(w, \sum_\alpha \varphi_\alpha T_* J_\alpha w) = \sum_\alpha \varphi_\alpha \omega_F(w, T_* J_\alpha w) > 0,
\ee
which is a contradiction. We conclude that any real eigenvector of $A$ must lie in $\ker T$. Let $x \in X$ be given. As $T$-regular points are always $T$-wrapped, and by hypothesis the same holds for all $T$-critical points, we know that $x$ is $T$-wrapped. Recall the subspace $P \subset \ker T \subset E$ used in \autoref{defn:wrapped}. We construct a decomposition ${\rm span}_\R P_x = \bigoplus_j \Pi_j$, with each $\Pi_j$ a real two-plane which is a $J_\alpha$-complex line for all $J_\alpha$ defined on $E_x$. Let $v \in P_x$ be given. By definition of $P$, there exists a sequence $(x_i)_{i \in \N}$ of $T$-regular points converging to $x$ and elements $v_i \in \ker T_{x_i}$ such that $v = \lim_{i\to\infty} v_i$. As the points $x_i$ are $T$-regular, the subspaces $\ker T_{x_i}$ are two-planes in $E_{x_i}$ oriented by the fiber-first convention. Pass to a subsequence so that the $\ker T_{x_i}$ converge to an oriented two-plane $\Pi \subset P_x$ containing $v$. Consider an open $W_\alpha$ containing $x$. Then there exists an $N_\alpha \in \N$ such that $x_i \in W_\alpha$ for all $i \geq N_\alpha$. But then for all $i \geq N_\alpha$, $\ker T_{x_i}$ is a $J_\alpha$-complex line, hence so is their limit $\Pi$. We conclude that $\Pi$ is a $J_\alpha$-complex line for each $J_\alpha$ defined at $x$. Proceed by induction to constructed multiple such real oriented two-planes $\Pi_j \subset P_x$ so that ${\rm span}_\R P_x = \bigoplus_j \Pi_j$, with each $\Pi_j$ being a $J_\alpha$-complex line for all $J_\alpha$ defined at $x$. Consider the quotient $Q_x := \ker T_x / {\rm span}_\R P_x$, which inherits an orientation from $\ker T_x$, which in turn is oriented as it is $J_\alpha$-complex for all $\alpha$ defined at $x$, all of which are $(\omega_F,T)$-tame. Then $Q_x$ inherits complex structures $\ol{J}_\alpha$ from each $J_\alpha$ defined at $x$, and these are all compatible with the orientation on $Q_x$. As $x$ is $T$-wrapped, $\dim_\C Q_x \leq 1$. But then there exists a fixed nondegenerate skew-symmetric bilinear form $\omega_x$ on $Q_x$ so that all $\ol{J}_\alpha$ are $\omega_x$-tame, as one can just pick an $\omega_x$ giving the orientation on $Q_x$. Consider the map $\ol{A}_x := \sum_\alpha \varphi_\alpha(x) \ol{J}_\alpha:\, Q_x \to Q_x$. Then $\ol{A}_x$ has no real eigenvalues on $Q_x$. As before, if $\ol{A}_x \ol{v} = \lambda \ol{v}$ for $0 \neq \ol{v} \in Q_x$ with $\lambda$ real, we would have
\be
	0 = \omega_x(\ol{v},\lambda \ol{v}) = \omega_x(\ol{v}, \ol{A}_x \ol{v}) = \omega_x(\ol{v}, \sum\limits_\alpha \varphi_\alpha(x) \ol{J}_\alpha \ol{v}) = \sum\limits_\alpha \varphi_\alpha(x) \omega_x(\ol{v}, \ol{J}_\alpha \ol{v}) > 0.
\ee
As $\ker T_x \cong {\rm span}_\R P_x \oplus Q_x$, we conclude that any real eigenvector of $A_x$ must lie in ${\rm span}_\R P_x = \bigoplus_j \Pi_j$. Construct a direct sum two-form $\wt{\omega}_x = \bigoplus_j \omega_j$ on ${\rm span}_\R P_x$ which tames each $J_\alpha$ at $x$. Then if $A_x v = \lambda v$ for $0 \neq v \in {\rm span}_\R P_x$ with $\lambda$ real and $v = \bigoplus_j v_j$ with respect to the direct sum decomposition of ${\rm span}_\R P_x$,
\be
	0 = \wt{\omega}_x(v, \lambda v) = \wt{\omega}_x(v, A_x v) = \sum\limits_j \omega_j(v_j, A_x v_j) = \sum\limits_{j,\alpha} \varphi_\alpha(x) \omega_j(v_j, J_\alpha(x) v_j) > 0.
\ee
We conclude that $A_x$ has no real eigenvalues, hence nor does $A$. By \autoref{lem:acsretraction} we obtain from $A$ a continuous complex structure $J = j(A)$ on $E$. As $T j(A) = j(B) T$ and convex combinations of $(\omega_F, T)$-tame endomorphisms are still tamed, the resulting \acs $J$ is $(\omega_F, T)$-tame.
\ep
\bp[ of \autoref{thm:bhyperfibration}] By \autoref{prop:blogbsymp} we need to show the existence of a $b$-symplectic form on $(X,Z_X)$. This is done by appealing to \autoref{thm:bthurstontrick}, hence it suffices to construct local closed two-forms $\eta_y$ around fibers so that the respective $b$-forms $\rho_X^* \eta_y$ tame a global \bacs $\btwo{J}$ on $\ker \b{df}$, and such that they are all cohomologous to the restriction of one global class $c \in H^2(X;\R)$. Construct around each point $x \in X$ a neighbourhood $W_x$ and a \bacs $\btwo{J}_x$ on $\b{TX}|_{W_x}$ which is $(\btwo{\omega}_Y, f)$-tame. These exist by definition around critical points of $X$ and away from critical points these are constructed as in the proof of \autoref{thm:bthurstonfibration}, using that the fibers of $f$ are two-dimensional. Apply \autoref{prop:acsconstruct} to the situation where $E = \b{TX}$, $F = f^*\b{TY}$, $T = \b{df}$ and $\omega_F = f^* \btwo{\omega}_Y$, to obtain a global \bacs $\btwo{J}$ on $X$ which is $(\btwo{\omega}_Y, f)$-tame. Let $S$ be the finite set of sections in the hypothesis, and define $c$ to be the Poincar\'e dual of $\sum_{s \in S} [{\rm im}(s)] \in H_{2n-2}(X;\R)$. To apply Poincar\'e duality in the absence of an orientation on $X$, the images of the sections must be cooriented. However, this exactly means that the fibers of $f$ must be oriented. Regular fibers obtain an orientation through $\btwo{J}$ and \autoref{prop:kerbdf}, while singular fibers are oriented in their smooth locus.

Let $y \in Y$ be given. By the definition of a hyperfibration, $\Delta_y$ is finite. Let $C_y \subset X$ be the disjoint union of closed balls centered at each point in $\Delta_y$. Choose a closed two-form $\sigma \in \Omega^2(C_y)$ so that $\rho_X^* \sigma$ tames $\btwo{J}$ on $\ker \b{df}|_{\Delta_y}$, noting that this is a condition at a finite set of points hence can easily be satisfied. Then, as taming $\btwo{J}$ on $\ker \b{df}$ is an open condition by \autoref{rem:tamingopen}, after possibly shrinking the balls in $C_y$ we can assume that $\btwo{J}|_{\ker \b{df}}$ is $\rho_X^* \sigma$-tame on the entirety of $C_y$. As by assumption $\btwo{J}$ is $(\omega_Y,f)$-tame, $F_y \setminus \Delta_y$ is a smooth noncompact $\btwo{J}$-holomorphic curve in $X \setminus \Delta_y$, whose complex $b$-orientation from $\btwo{J}$ agrees with its preimage $b$-orientation. By assumption there exists for each component of $F_y \setminus \Delta_y$ a section in $S$ intersecting that component positively. Choose $C_y$ so that $\partial C_y$ is transverse to $F_y$ and consider the intersecting circles $F_y \cap \partial C_y$. Connect each such circle to these points of intersection by a path in $F_y \setminus \Delta_y$. Let $C_y^0$ be the disjoint union of smaller concentric closed balls around $\Delta_y$ disjoint from these paths and with $\partial C_y^0$ transverse to $F_y$. Then each component $F_i$ of the compact surface $F_y \setminus {\rm int}(C_y^0)$ either lies fully inside ${\rm int}(C_y)$, or has a point of intersection with a section in $S$. In the latter case we will say that $F_i$ \emph{intersects $S$}. See the following figures for illustration.

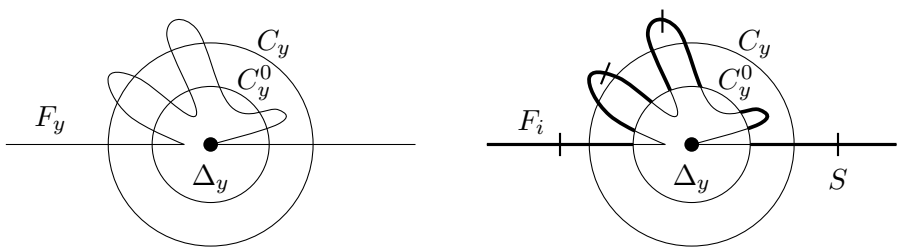
\newcommand{\fibercurve}[2][black]{
	\draw[line width = #2, #1]
	plot[smooth, tension=.7] coordinates{($(a) + (-0.9,0)$) ($(a) + (-2.8,1)$) ($(a) + (-3.5,2.3)$) ($(a) + (-2.5,2.3)$) ($(a) + (-0.5,1)$) ($(a) + (-1.5,3.8)$) ($(a) + (-.5,4)$) ($(a) + (0.7,1.3)$) ($(a) + (2.3,1.2)$) ($(a) + (2.3,.7)$) (a)};
}
\begin{center}
\begin{tikzpicture}[scale=1,x=1em,y=1em]
\begin{scope}[xshift=-180pt]
	\coordinate (a) at (0,0);
	\fill (a) circle (0.25) node [label=below:{$\Delta_y$}]{};
	\draw[thin]
		(a) circle (2) node [label={[label distance = 0.9em]65:{$C_y^0$}}]{}
		(a) circle (3.5) node [label={[label distance = 2.5em]65:{$C_y$}}]{}
		($(a) + (-7,0)$) -- ($(a) + (-0.9,0)$) node [above, near start] {$F_y$}
		(a) -- ($(a) + (7,0)$);
	\fibercurve{0.4pt}
\end{scope}
\begin{scope}
	\coordinate (a) at (0,0);
	\fill (a) circle (0.25) node [label=below:{$\Delta_y$}]{};
	\draw[thin]
		(a) circle (2) node [label={[label distance = 0.9em]65:{$C_y^0$}}]{}
		(a) circle (3.5) node [label={[label distance = 2.5em]65:{$C_y$}}]{}
		($(a) + (-7,0)$) -- ($(a) + (-0.9,0)$) node [above, near start] {$F_i$}
		(a) -- ($(a) + (7,0)$);
	\fibercurve{0.4pt}
	\draw[very thick]
		($(a) + (-7,0)$) -- ($(a) + (-2,0)$)
		($(a) + (2,0)$) -- ($(a) + (7,0)$);
	\draw[thick]
		($(a) + (-3.1,2.075)$) -- ($(a) + (-2.8,2.825)$) 
		($(a) + (-1,3.85)$) -- ($(a) + (-1,4.65)$) 
		($(a) + (-4.5,-0.4)$) -- ($(a) + (-4.5,0.4)$) 
		($(a) + (+5,-0.4)$) -- ($(a) + (+5,0.4)$) node [label={[label distance = 0.5em]below:{$S$}}]{}; 
\end{scope}
\begin{scope}[even odd rule]
	\coordinate (a) at (0,0);
	\clip ($(a) + (-7,-3.5)$) rectangle ($(a) + (7,5)$)
	(a) circle (2);
	\fibercurve{1.4pt}
\end{scope}
\end{tikzpicture}
\end{center}

Let $W_y$ be the union of ${\rm int}(C_y^0)$ with a tubular neighbourhood rel boundary of $F_y \setminus {\rm int}(C_y^0)$ inside $X \setminus {\rm int}(C_y^0)$. Extend each $F_i$ to a closed oriented smooth surface $\wh{F}_i \subset W_y$ by arbitrarily attaching a surface inside $C_y^0$. Then the classes $[\wh{F}_i] \in H_2(W_y;\Z)$ form a basis for the homology of $W_y$. Indeed, contracting the whole neighbourhood $C_y^0$ to $\Delta_y$ we see that $F_y$ becomes homotopy equivalent to a wedge of the $F_i$. This $W_y$ is the desired neighbourhood of $F_y$ on which we construct the form $\eta_y$ in the hypothesis of \autoref{thm:bthurstontrick}, using ideas similar to the proof of \autoref{thm:blfblog}. Since $F_y$ is $\btwo{J}$-holomorphic with $\btwo{J}|_{\ker \b{df}}$ being $\rho_X^*\sigma$-tame on $C_y$, $\sigma|_{F_i \cap C_y}$ is a positive area form for each $i$. For each $F_i$ intersecting $S$, let $\sigma_i$ be an extension of $\sigma$ to $F_i$ as a positive area form with total area $\langle \sigma_i, [\wh{F}_i] \rangle$ equal to $\#(F_i \cap \im(S)) > 0$. Let $p$ denote the tubular neighbourhood projection onto $F_y \setminus {\rm int}(C_y^0)$ and let $f:\, C_y \to [0,1]$ be a smooth radial function defined on each ball around points in $\Delta_y$ so that $f \equiv 0$ in a neighbourhood of $C_y^0$ and $f \equiv 1$ in a neighbourhood of $\partial C_y$, which is smoothly extended to $W_y$ by being identically $1$ outside of $C_y$. See the following figure for illustration.

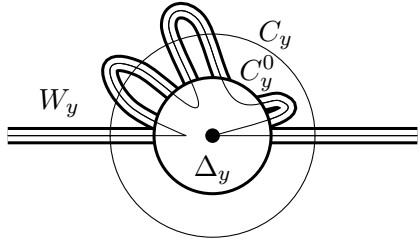
\begin{center}
\begin{tikzpicture}[scale=1,x=1em,y=1em]
\coordinate (a) at (0,0);

\begin{scope}[even odd rule]
\clip ($(a) + (-7.1,-3.5)$) rectangle ($(a) + (7.1,5)$)
(a) circle (2);
\draw[line width = 2.3mm]
	($(a) + (-7,0)$) -- ($(a) + (7,0)$) node [above, very near start] {$W_y$};
\draw[line width = 1.5mm, white]
	($(a) + (-7,0)$) -- ($(a) + (7,0)$);
\draw[thin]
	($(a) + (-7,0)$) -- ($(a) + (7,0)$);
\end{scope}

\begin{scope}[even odd rule]
\coordinate (a) at (0,0);
\clip ($(a) + (-7,-3.5)$) rectangle ($(a) + (7,5)$)
(a) circle (2);

\fibercurve{2.3mm}
\fibercurve[white]{1.5mm}
\fibercurve{0.4pt}
\end{scope}

\begin{scope}
\fill (a) circle (0.25) node [label=below:{$\Delta_y$}]{};
\draw[line width=0.4mm]
	(a) circle (2) node [label={[label distance = 0.85em]65:{$C_y^0$}}]{};
\draw[thin]
	(a) circle (3.5) node [label={[label distance = 2.5em]65:{$C_y$}}]{}
	($(a) + (-7,0)$) -- ($(a) + (-0.9,0)$) 
	(a) -- ($(a) + (7,0)$);
\fibercurve{0.4pt}
\end{scope}
\end{tikzpicture}
\end{center}

On the balls $C_y$, the form $\sigma$ is exact, say equal to $\sigma = d\alpha$ for $\alpha \in \Omega^1(C_y)$. Define a two-form $\eta_y$ on $W_y$ by $\eta_y := \sum_i p^* (f \sigma_i) + d((1-f) \alpha)$. In other words, $\sigma$ is extended by $0$ outside $C_y$, while the $\sigma_i$ are extended by $0$ inside $C_y^0$. For all functions $f \neq 0$ the form $f \sigma_i$ is a closed area form defined on the surface $F_i$, so that $\eta_y$ is nonnegative and closed. On $W_y \cap C_y^0 = C_y^0$ we have $f \equiv 0$ so that $\eta_y = d\alpha = \sigma$, hence there $\rho_X^* \eta_y = \rho_X^* \sigma$ tames $\btwo{J}|_{\ker \b{df}}$. Similarly, on $F_y \setminus \Delta_y$ the forms $\sigma_i$ are area forms for the orientation given by $\btwo{J}$. But then $\rho_X^* \eta_y$ tames $\btwo{J}$ on $T_x F_y = \ker df_x \cong \ker \b{df}_x$ for all $x \in F_y \setminus \Delta_y$ as this condition is convex. By openness of the taming condition, narrowing the tubular neighbourhood $p$ defining $W_y$ ensures that $\rho_X^* \eta_y$ tames $\btwo{J}|_{\ker \b{df}_x}$ for all $x \in W_y \setminus {\rm int}(C_y^0)$. But then $\rho_X^* \eta_y$ tames $\btwo{J}$ on $\ker \b{df}_x$ for all $x \in W_y$.

What remains is to show that $c|_{W_y} = [\eta_y] \in H^2(W_y; \R)$. Recall that every component $F_i$ of $F_y \setminus {\rm int}(C_y^0)$ either intersects $S$ or lies in ${\rm int}(C_y)$. For those $F_i$ intersecting $S$ we have $\langle [\eta_y], [\wh{F}_i] \rangle = \langle [\sigma_i], [\wh{F}_i] \rangle = \# (F_i \cap {\rm im}(S)) = \langle c, [\wh{F}_i] \rangle$. For $F_i$ disjoint from $S$ we know that $F_i \subset C_y$, but $\eta_y$ is exact in $C_y$, so that we have $\langle [\eta_y], [\wh{F}_i] \rangle = \langle 0, [\wh{F}_i] \rangle = 0 = \langle c, [\wh{F_i}] \rangle$. As $[\eta_y]$ agrees with the class $c|_{W_y} \in H^2(W_y; \R)$ when evaluated on the classes [$\wh{F}_i$], which form a basis of $H_2(W_y;\Z)$, we see that $[\eta_y] = c|_{W_y} \in H^2(W_y; \R)$ as desired. Applying \autoref{thm:bthurstontrick} we conclude that $(X,Z_X)$ admits a \blog.
\ep
%
%
%
\section{Log-symplectic and \fsymp{s}}
\label{sec:blogsandfsymps}
In this section we discuss \fsymp{}s and show that \blog{}s naturally give rise to \fsymp{}s. \Fsymp{}s are studied amongst others in \cite{Baykur06, CannasDaSilvaGuilleminPires11, CannasDaSilva10, CannasDaSilvaGuilleminWoodward00}.

\begin{defn} A \emph{\fsymp} on a compact $2n$-dimensional manifold $X$ is a closed two-form $\omega$ such that $\bigwedge^n \omega$ is transverse to the zero section in $\bigwedge^{2n} T^*X$, and such that $\omega^{n-1}|_{Z_\omega} \neq 0$, where $Z_\omega = (\bigwedge^n \omega)^{-1}(0)$. The hypersurface $Z_\omega$ is called the \emph{folding locus} of $\omega$, while its complement $X \setminus Z_\omega$ is called the \emph{symplectic locus}. A \fsymp is called \emph{bona fide} if $Z_\omega \neq \emptyset$.
\end{defn}
This definition should be compared with \autoref{defn:blog}. We will say the pair $(X,Z_X)$ \emph{admits a \fsymp} if $X$ admits a \fsymp $\omega$ for which $Z_\omega = Z_X$. According to the Darboux model, a \fsymp $\omega$ is locally given by $\omega = x_1 d x_1 \wedge dx_2 + \dots + dx_{2n-1} \wedge dx_{2n}$, using coordinates $x_i$ in a neighbourhood $U$ such that $Z_\omega \cap U = \{x_1 = 0\}$.

Cannas da Silva gave a homotopical characterization for an orientable manifold to admit a \fsymp.
\begin{thm}[\cite{CannasDaSilva10}]\label{thm:fsympstableacs} Let $X$ be an orientable manifold. Then $X$ admits a \fsymp if and only if $X$ admits a stable \acs. In particular, every orientable four-manifold admits a \fsymp.
\end{thm}
Further, Baykur has given a construction showing there is a relation between \alf{}s on four-manifolds and \fsymp{}s. This should be compared with \autoref{thm:alflogsymp}.
\begin{thm}[{\cite[Proposition 3.2]{Baykur06}}]\label{thm:alffsymp} Let $f:\, X^4 \to \Sigma^2$ be an \alf between compact connected manifolds. Assume that the generic fiber $F$ is orientable and $[F] \neq 0 \in H_2(X;\R)$. Then $X$ admits a \fsymp.
\end{thm}
From our point of view this theorem does not come as a surprise. Indeed, every \blog gives rise to a \fsymp, as is known in the Poisson community. We learned the proof of this result from M{\u{a}}rcu{\c{t}} and Frejlich.
\begin{thm}\label{thm:blogfolded} Let $(X^{2n},Z_X,\pi)$ be a compact log-symplectic manifold. Then $(X,Z_X)$ admits a \fsymp $\omega$ for which $\omega = \pi^{-1}$ outside a neighbourhood of $Z_X$.
\end{thm}
\bp By \autoref{prop:bloglocalform} a neighbourhood of each connected component $Z$ of $Z_X$ is equivalent to a neighbourhood $U$ of the zero section of the normal bundle $NZ$ equipped with a distance function $|x|$, so that $\pi^{-1} = d \log |x| \wedge \theta + \sigma$ for closed one- and two-forms $\theta$ and $\sigma$ on $Z$ satisfying $\theta \wedge \sigma^{n-1} \neq 0$. By rescaling $|x|$ we can assume that $U$ contains all points of distance at most $e^2 + 1$ away from the zero section. Denote $\omega_Z = d |x|^2 \wedge \theta + \sigma$ and let $f:\, \R_+ \to \R$ be a smooth monotone interpolation between the functions $f_0:\, [0,1] \to \R$, $f_0(x) = x^2$ and $f_1:\, [e^2,\infty) \to \R$, $f_1(x) = \log x$. Consider the closed two-form $\omega_f = df(|x|) \wedge \theta + \sigma$, extended by $\pi^{-1}$ outside of $U$. Then $\omega_f = \pi^{-1}$ away from $Z$, while near $Z$ we have $\omega_f = \omega_Z$. Moreover, $\omega_f$ is symplectic on $X \setminus Z$ by monotonicity of $f$. Perform this procedure for all connected components of $Z_X$ to obtain a closed two-form $\omega$ on $X$ for which $\omega = \pi^{-1}$ away from $Z_X$. By the local description near $Z_X$ it follows that $\omega^n$ vanishes transversally with $Z_\omega = Z_X$. Further, the restriction of $\omega^{n-1}$ to $Z_\omega$ is equal to $\sigma^{n-1}$, hence is nonvanishing as $\theta \wedge \sigma^{n-1} \neq 0$. We conclude that $\omega$ is a \fsymp for $(X,Z_X)$.
\ep
The previous theorem, together with \autoref{thm:alflogsymp}, implies \autoref{thm:alffsymp}. Moreover, the \fsymp of \autoref{thm:alffsymp} agrees with the one obtained through our methods, as is hinted at by the fact that in Baykur's construction the folding locus fibers over the circle.

The converse to \autoref{thm:blogfolded} does not hold. For example, $S^4$ does not admit a \blog by \autoref{thm:blogobstr}, while it does admit \fsymp{}s by \autoref{thm:fsympstableacs}. Similarly, by \autoref{thm:bonafideobstr} and results from Seiberg-Witten theory due to Taubes \cite{Taubes94}, the four-manifold $\C P^2 \# \C P^2$ does not admit \blog{}s, bona fide or not. However, by \autoref{thm:fsympstableacs} it does admit a \fsymp.

The reason the converse to \autoref{thm:blogfolded} is false is essentially because the information contained in the one-form $\theta$ determined by the \blog (see \autoref{prop:bloglocalform}) is lost when passing to the folded-symplectic world. Note here that a \fsymp $\omega$ restricts to $Z_\omega$ to define a one-dimensional foliation $\ker(\left.\omega\right|_{Z_\omega})$ called the \emph{null foliation}, and $\left.\omega\right|_{Z_\omega}$ is a \emph{pre-symplectic structure} on $Z_\omega$, i.e.\ a closed two-form of maximal rank. On the other hand, a \blog $\pi$ induces a cosymplectic structure on $Z_\pi$ by \autoref{prop:bloglocalform}, so that the associated nowhere-vanishing closed one-form gives a codimension-one foliation on $Z_\pi$. Further, there is a symplectic structure on the leaves.

The following result makes precise that it is exactly the existence of a suitable closed one-form $\theta$ on $Z_X$ that makes the converse to \autoref{thm:blogfolded} hold.
\begin{thm}\label{thm:foldedblog} Let $(X^{2n},Z_X,\omega)$ be a compact folded-symplectic manifold. Assume that there exists a closed one-form $\theta \in \Omega^1(Z_X)$ such that $\theta \wedge \omega^{n-1}|_{Z_X} \neq 0$. Then $(X,Z_X)$ admits a \blog $\pi$ for which $\pi = \omega^{-1}$ outside a neighbourhood of $Z_X$.
\end{thm}
In other words, when the folded-symplectic form $\omega$ can be complemented to give a cosymplectic structure on $Z_X$, one can turn $\omega$ into a \blog.
\bp Consider the normal bundle $NZ_X$ and let $U$ be a neighbourhood of the zero section. Choose a distance function $|x|$ for $Z_X$ which is constant outside of $\overline{U}$. Note that $d \log |x| \wedge \theta \wedge \omega^{n-1}$ is nonzero at $Z_X$ as $\theta \wedge \omega^{n-1}|_{Z_X} \neq 0$ and $|x|$ is transverse to $Z_X$. By continuity it is still nonzero in a neighbourhood $V \subset U$ of $Z_X$. Let $f = f(|x|)$ be a smooth function on $NZ_X$ so that $f \equiv 1$ on $V$ and $f \equiv 0$ near $\partial U$, which is then extended to $X$ by being identically $0$ outside $U$. Define a closed $b$-two-form $\omega_f = t \, d (f \log |x|) \wedge \theta + \omega \in \b{\Omega}^2(X)$ for the $b$-manifold $(X,Z_X)$, where $t \neq 0$ is a real parameter. Choose the sign of $t$ so that the forms $t \, d(f \log|x|) \wedge \theta$ and $\omega$ give the same orientation on $V \setminus Z_X$. We have $\omega_f^n = t n \, d (f \log |x|) \wedge \theta \wedge \omega^{n-1} + \omega^n$, so by choosing $t$ small enough we conclude that $\omega_f$ is a $b$-symplectic form for which $\omega_f = \omega$ outside $\overline{U}$. By \autoref{prop:blogbsymp} the dual bivector $\pi$ to $\omega_f$ is a \blog for $(X,Z_X)$, and $\pi = \omega^{-1}$ away from $Z_X$.
\ep
%
%
\bibliographystyle{hyperamsplain-nodash}
\bibliography{blogsfibrations}

\providecommand{\bysame}{\leavevmode\hbox to3em{\hrulefill}\thinspace}
\providecommand{\MR}{\relax\ifhmode\unskip\space\fi MR }
\providecommand{\MRhref}[2]{%
  \href{http://www.ams.org/mathscinet-getitem?mr=#1}{#2}
}
\providecommand{\href}[2]{#2}
\begin{thebibliography}{10}

\bibitem{AkbulutKarakurt08}
S.~Akbulut and {\c{C}}.~Karakurt, \emph{Every 4-manifold is {BLF}},
  \href{http://www.gokovagt.org/journal/2008/jggt08-akbulutkarakurt.pdf}{J.
  G\"okova Geom. Topol. GGT} \textbf{2} (2008), 83--106.

\bibitem{AurouxDonaldsonKatzarkov05}
D.~Auroux, S.~K. Donaldson, and L.~Katzarkov, \emph{Singular {L}efschetz
  pencils}, \href{http://dx.doi.org/10.2140/gt.2005.9.1043}{Geom. Topol.
  \textbf{9} (2005)}, 1043--1114.

\bibitem{Baykur06}
R.~{\.I}. Baykur, \emph{K\"ahler decomposition of 4-manifolds},
  \href{http://dx.doi.org/10.2140/agt.2006.6.1239}{Algebr. Geom. Topol.
  \textbf{6} (2006)}, 1239--1265.

\bibitem{Baykur08}
R.~{\.I}. Baykur, \emph{Existence of broken {L}efschetz fibrations},
  \href{http://dx.doi.org/10.1093/imrn/rnn101}{Int. Math. Res. Not. IMRN
  (2008)}, Art. ID rnn 101, 15.

\bibitem{Baykur09}
R.~{\.I}. Baykur, \emph{Topology of broken {L}efschetz fibrations and
  near-symplectic four-manifolds},
  \href{http://dx.doi.org/10.2140/pjm.2009.240.201}{Pacific J. Math.
  \textbf{240} (2009)}, no.~2, 201--230.

\bibitem{CannasDaSilvaGuilleminPires11}
A.~Cannas~da Silva, V.~Guillemin, and A.~R. Pires, \emph{Symplectic origami},
  \href{http://dx.doi.org/10.1093/imrn/rnq241}{Int. Math. Res. Not. IMRN
  (2011)}, no.~18, 4252--4293.

\bibitem{CannasDaSilva10}
A.~Cannas~da Silva, \emph{Fold-forms for four-folds},
  \href{http://dx.doi.org/10.4310/JSG.2010.v8.n2.a3}{J. Symplectic Geom.
  \textbf{8} (2010)}, no.~2, 189--203.

\bibitem{CannasDaSilvaGuilleminWoodward00}
A.~Cannas~da Silva, V.~Guillemin, and C.~Woodward, \emph{On the unfolding of
  folded symplectic structures},
  \href{http://dx.doi.org/10.4310/MRL.2000.v7.n1.a4}{Math. Res. Lett.
  \textbf{7} (2000)}, no.~1, 35--53.

\bibitem{Cavalcanti13}
G.~R. {Cavalcanti}, \emph{{Examples and counter-examples of log-symplectic
  manifolds}}, \href{http://arxiv.org/abs/1303.6420}{{\tt arXiv:1303.6420}}.

\bibitem{Donaldson99}
S.~K. Donaldson, \emph{Lefschetz pencils on symplectic manifolds},
  \href{http://projecteuclid.org/euclid.jdg/1214425535}{J. Differential Geom.}
  \textbf{53} (1999), no.~2, 205--236.

\bibitem{EtnyreFuller06}
J.~B. Etnyre and T.~Fuller, \emph{Realizing 4-manifolds as achiral {L}efschetz
  fibrations}, \href{http://dx.doi.org/10.1155/IMRN/2006/70272}{Int. Math. Res.
  Not. (2006)}, Art. ID 70272, 21.

\bibitem{FrejlichMartinezTorresMiranda15}
P.~{Frejlich}, D.~{Mart{\'{\i}}nez Torres}, and E.~{Miranda}, \emph{{A note on
  symplectic topology of $b$-manifolds}},
  \href{http://arxiv.org/abs/1312.7329}{{\tt arXiv:1312.7329}}.

\bibitem{GayKirby07}
D.~T. Gay and R.~Kirby, \emph{Constructing {L}efschetz-type fibrations on
  four-manifolds}, \href{http://dx.doi.org/10.2140/gt.2007.11.2075}{Geom.
  Topol. \textbf{11} (2007)}, 2075--2115.

\bibitem{Gompf94}
R.~E. Gompf, \emph{A new construction of symplectic manifolds},
  \href{http://dx.doi.org/10.2307/2118554}{Ann. of Math. (2) \textbf{142}
  (1995)}, no.~3, 527--595.

\bibitem{Gompf01}
R.~E. Gompf, \emph{The topology of symplectic manifolds},
  \href{http://journals.tubitak.gov.tr/math/issues/mat-01-25-1/mat-25-1-2-0103-2.pdf}{Turkish
  J. Math.} \textbf{25} (2001), no.~1, 43--59.

\bibitem{Gompf04two}
R.~E. Gompf, \href{http://dx.doi.org/10.2140/gtm.2004.7.267}{\emph{Symplectic
  structures from {L}efschetz pencils in high dimensions}}, Proceedings of the
  {C}asson {F}est, Geom. Topol. Monogr., vol.~7, Geom. Topol. Publ., Coventry,
  2004, pp.~267--290 (electronic).

\bibitem{Gompf04}
R.~E. Gompf, \emph{Toward a topological characterization of symplectic
  manifolds}, \href{http://dx.doi.org/10.4310/JSG.2004.v2.n2.a1}{J. Symplectic
  Geom. \textbf{2} (2004)}, no.~2, 177--206.

\bibitem{Gompf05}
R.~E. Gompf, \emph{Locally holomorphic maps yield symplectic structures},
  \href{http://dx.doi.org/10.4310/CAG.2005.v13.n3.a2}{Comm. Anal. Geom.
  \textbf{13} (2005)}, no.~3, 511--525.

\bibitem{GompfStipsicz99}
R.~E. Gompf and A.~I. Stipsicz,
  \href{http://dx.doi.org/10.1090/gsm/020}{\emph{{$4$}-manifolds and {K}irby
  calculus}}, Graduate Studies in Mathematics, vol.~20, American Mathematical
  Society, Providence, RI, 1999.

\bibitem{GMP13}
V.~Guillemin, E.~Miranda, and A.~R. Pires, \emph{Symplectic and {P}oisson
  geometry on {$b$}-manifolds},
  \href{http://dx.doi.org/10.1016/j.aim.2014.07.032}{Adv. Math. \textbf{264}
  (2014)}, 864--896.

\bibitem{Honda04}
K.~Honda, \emph{Local properties of self-dual harmonic 2-forms on a
  4-manifold}, \href{http://dx.doi.org/10.1515/crll.2004.2004.577.105}{J. Reine
  Angew. Math. \textbf{577} (2004)}, 105--116.

\bibitem{Klaasse16}
R.~L. Klaasse, Ph.D. thesis, Utrecht University, 2016. Work in progress.

\bibitem{LeBrun97}
C.~LeBrun, \emph{Yamabe constants and the perturbed {S}eiberg-{W}itten
  equations}, \href{http://arxiv.org/abs/dg-ga/9605009}{Comm. Anal. Geom.}
  \textbf{5} (1997), no.~3, 535--553.

\bibitem{Lekili09}
Y.~Lekili, \emph{Wrinkled fibrations on near-symplectic manifolds},
  \href{http://dx.doi.org/10.2140/gt.2009.13.277}{Geom. Topol. \textbf{13}
  (2009)}, no.~1, 277--318. Appendix B by R. {\.I}nan{\c{c}} Baykur.

\bibitem{MarcutOsorno14}
I.~M{\u{a}}rcu{\c{t}} and B.~Osorno~Torres, \emph{Deformations of
  log-symplectic structures}, \href{http://dx.doi.org/10.1112/jlms/jdu023}{J.
  Lond. Math. Soc. (2) \textbf{90} (2014)}, no.~1, 197--212.

\bibitem{MarcutOsorno14two}
I.~M{\u{a}}rcu{\c{t}} and B.~Osorno~Torres, \emph{On cohomological obstructions
  for the existence of log-symplectic structures},
  \href{http://dx.doi.org/10.4310/JSG.2014.v12.n4.a6}{J. Symplectic Geom.
  \textbf{12} (2014)}, no.~4, 863--866.

\bibitem{McDuffSalamon98}
D.~McDuff and D.~Salamon,
  \emph{\href{https://global.oup.com/academic/product/introduction-to-symplectic-topology-9780198504511}{Introduction
  to symplectic topology}}, second ed., Oxford Mathematical Monographs, The
  Clarendon Press, Oxford University Press, New York, 1998.

\bibitem{Melrose93}
R.~B. Melrose,
  \emph{\href{https://www.crcpress.com/The-Atiyah-Patodi-Singer-Index-Theorem/Melrose/9781568810027}{The
  {A}tiyah-{P}atodi-{S}inger index theorem}}, Research Notes in Mathematics,
  vol.~4, A K Peters, Ltd., Wellesley, MA, 1993.

\bibitem{Radko02}
O.~Radko, \emph{A classification of topologically stable {P}oisson structures
  on a compact oriented surface},
  \href{http://projecteuclid.org/euclid.jsg/1092403031}{J. Symplectic Geom.}
  \textbf{1} (2002), no.~3, 523--542.

\bibitem{Taubes94}
C.~H. Taubes, \emph{The {S}eiberg-{W}itten invariants and symplectic forms},
  \href{http://dx.doi.org/10.4310/MRL.1994.v1.n6.a15}{Math. Res. Lett.
  \textbf{1} (1994)}, no.~6, 809--822.

\bibitem{Thurston76}
W.~P. Thurston, \emph{Some simple examples of symplectic manifolds},
  \href{http://dx.doi.org/10.2307/2041749}{Proc. Amer. Math. Soc. \textbf{55}
  (1976)}, no.~2, 467--468.

\end{thebibliography}
\end{document}